# VERSAL EMBEDDINGS OF COMPACT, REGULAR 3-PSEUDOCONCAVE CR SUBMANIFOLDS

PETER L. POLYAKOV

ABSTRACT. We prove that for an induced CR structure on a compact, generic, regular 3-pseudoconcave CR submanifold $\mathbf{M} \subset \mathbf{G}$, of a complex manifold $\mathbf{G}$, satisfying condition $\dim H^1\left(\mathbf{M}, T'(\mathbf{G})|_{\mathbf{M}}\right) = 0$ all the close CR structures are induced by close embeddings.

## 1. INTRODUCTION.

Let $\mathbf{M}$ be a generic CR submanifold in a complex manifold $\mathbf{G}$, i.e. such a submanifold, that for any $z \in \mathbf{M}$ there exist a neighborhood $\mathcal{U} \ni z$ in $\mathbf{G}$ and smooth real valued functions $\{\rho_k\}_1^m$ ($1 < m < n-1$) on $\mathcal{U}$ such that

$$\mathbf{M} \cap \mathcal{U} = \{z \in \mathbf{G} \cap \mathcal{U} : \rho_1(z) = \cdots = \rho_m(z) = 0\}, \qquad (1)$$
$$\partial \rho_1 \wedge \cdots \wedge \partial \rho_m \neq 0 \text{ on } \mathbf{M} \cap \mathcal{U}.$$

The induced CR structure $\mathbf{M}$ is defined by the subbundles

$$T''(\mathbf{M}) = T''(\mathbf{G})|_{\mathbf{M}} \cap \mathbf{C}T(\mathbf{M}) \quad \text{and} \quad T'(\mathbf{M}) = T'(\mathbf{G})|_{\mathbf{M}} \cap \mathbf{C}T(\mathbf{M}),$$

where $\mathbf{C}T(\mathbf{M})$ is the complexified tangent bundle of $\mathbf{M}$ and the subbundles $T''(\mathbf{G})$ and $T'(\mathbf{G}) = \overline{T''(\mathbf{G})}$ of the complexified tangent bundle $\mathbf{C}T(\mathbf{G})$ define a complex structure on $\mathbf{G}$.

If we fix a Hermitian scalar product on $\mathbf{G}$ then there exists a subbundle $\mathbf{N} \in \mathbf{C}T(\mathbf{M})$ of complex dimension $m$ such that

$$\mathbf{C}T(\mathbf{M}) = T'(\mathbf{M}) \oplus T''(\mathbf{M}) \oplus \mathbf{N} \text{ with } T'(\mathbf{M}) \perp \mathbf{N} \text{ and } T''(\mathbf{M}) \perp \mathbf{N}.$$

The Levi form of $\mathbf{M}$ is defined (cf. [He2]) as a Hermitian form on $T'(\mathbf{M})$ with values in $\mathbf{N}$

$$\mathcal{L}_z(L(z)) = \sqrt{-1} \cdot \pi\left(\left[\overline{L}, L\right]\right)(z) \quad \left(L(z) \in T'(\mathbf{M})_z\right),$$

where $\left[\overline{L}, L\right] = \overline{L}L - L\overline{L}$ and $\pi$ is the projection of $\mathbf{C}T(\mathbf{M})$ along $T'(\mathbf{M}) \oplus T''(\mathbf{M})$ onto $\mathbf{N}$.

The Levi form of $\mathbf{M}$ at the point $z \in \mathbf{M}$ in the direction of the unit vector $\theta = (\theta_1, \ldots \theta_m) \in \text{Re}\mathbf{N}_z$ is a scalar Hermitian form on $T'(\mathbf{M})_z$

$$\mathcal{L}_z^\theta(L(z)) = \langle \theta, \mathcal{L}_z(L) \rangle_z,$$

where $\langle \cdot, \cdot \rangle_z$ is a Hermitian scalar product.

Following [He2] we call $\mathbf{M}$ q-pseudoconcave (weakly q-pseudoconcave) at $z \in \mathbf{M}$ in the direction $\theta$ if the Levi form of $\mathbf{M}$ at $z$ in this direction has at least $q$ negative ($q$ nonpositive) eigenvalues and we call $\mathbf{M}$ q-pseudoconcave (weakly q-pseudoconcave) at $z \in \mathbf{M}$ if it is q-pseudoconcave (weakly q-pseudoconcave) in all directions.

We call a q-pseudoconcave CR manifold $\mathbf{M}$ a *regular q-pseudoconcave CR manifold* [P2] if for any $z \in \mathbf{M}$ there exist an open neighborhood $U \ni z$ in $\mathbf{M}$ and a family $E_q(\theta, z)$ of $q$-dimensional complex linear subspaces in $T_z^c(\mathbf{M})$ smoothly depending on $(\theta, z) \in \mathbb{S}^{m-1} \times U$ and such that the Levi form $\langle \theta, L_z(\mathbf{M}) \rangle$ is strictly negative on $E_q(\theta, z)$.

In this paper we consider small deformations of the induced CR structure on a compact,

generic, regular 3-pseudoconcave submanifold $\mathbf{M} \in \mathbf{G}$ of codimension $m$. We fix an embedding $\mathcal{E}_0 : \mathbb{M} \to \mathbf{M}_0 \subset \mathbf{G}$ of a compact $C^p$ manifold $\mathbb{M}$ into a complex manifold $\mathbf{G}$ such that $\mathbf{M}_0 = \mathcal{E}_0(\mathbb{M})$ is a generic, regular 3-pseudoconcave CR submanifold of $\mathbf{G}$. Let $\{\mathcal{U}^i\}_1^N$ be a finite cover of some neighborhood of $\mathbf{M}_0$ in $\mathbf{G}$ such that in each $\mathcal{U}^i$ the manifold $\mathbf{M}_0 \cap \mathcal{U}^i$ has the form (1) with defining functions $\left\{\rho_{i,l}^{(0)}\right\}_{1 \leq l \leq m}$. If $\mathcal{E}$ is an embedding of $\mathbb{M}$ into $\mathbf{G}$ close to $\mathcal{E}_0$ with $\mathbf{M} = \mathcal{E}(\mathbb{M})$ then the map $\mathcal{F} = \mathcal{E} \circ \mathcal{E}_0^{-1} : \mathbf{M}_0 \to \mathbf{M}$ may be defined in some small enough neighborhood $U^i = \mathbf{M}_0 \cap \mathcal{U}^i$ as $\mathcal{F}(z) = z + f_i(z)$, with $f_i \in [C^p(U^i)]^n$. If $|f|_p = \max_i\{|f_i|_p\}$ is small enough then $\mathcal{G} = \mathcal{F}^{-1} : \mathbf{M} \to \mathbf{M}_0$ is also well defined and has the form $\mathcal{G}(z) = z + g_i(z)$ in some neighborhood $V^i \subset \mathcal{F}(U^i)$ with $g_i \in [C^p(V^i)]^n$. We denote

$$|\mathcal{E}|_p = \max_{1 \leq i \leq N} \{|g_i|_p\}.$$

For a compact $C^p$ submanifold $\mathbf{M} \subset \mathbf{G}$ and a holomorphic vector bundle $\mathcal{B}$ on $\mathbf{G}$ we say that $\dim H^r\left(\mathbf{M}, \mathcal{B}|_\mathbf{M}\right) = 0$ for $r \in \mathbb{Z}^+$ if for any $\bar{\partial}_\mathbf{M}$-closed form $f \in C_{(0,r)}^k\left(\mathbf{M}, \mathcal{B}|_\mathbf{M}\right)$, $(1 \leq k \leq p)$, there exists a form $h \in C_{(0,r-1)}^k\left(\mathbf{M}, \mathcal{B}|_\mathbf{M}\right)$ such that $\bar{\partial}_\mathbf{M} h = f$.

In the theorem below we formulate the main result of the paper - versality of the family of CR structures induced by close embeddings for a fixed compact, generic, regular 3-pseudoconcave submanifold $\mathbf{M} \in \mathbf{G}$, satisfying condition $\dim H^1\left(\mathbf{M}, T'(\mathbf{G})|_\mathbf{M}\right) = 0$.

**Theorem 1.** *Let $k \in \mathbb{Z}$ be such that $k \geq 17$ and let $\mathcal{E} : \mathbb{M} \to \mathbf{G}$ be a $C^{2k+5}$-embedding of a compact manifold $\mathbb{M}$ into a complex manifold $\mathbf{G}$ such that $\mathbf{M} = \mathcal{E}(\mathbb{M})$ is a generic, regular 3-pseudoconcave CR manifold with $\dim H^1\left(\mathbf{M}, T'(\mathbf{G})|_\mathbf{M}\right) = 0$. Then any small $C^{2k+5}$-deformation of the induced CR structure $\mathbf{M}$ on $\mathbb{M}$ may be obtained as an induced CR structure from a $C^k$-embedding $\mathcal{F} : \mathbb{M} \to \mathbf{G}$, close to $\mathcal{E}$.*

The problem of embeddability for deformed CR structures on an embedded, compact, strictly pseudoconvex hypersurface was first considered by L. Boutet de Monvel in [B], where the preservation of local embeddability for such CR structures was proved. We refer to the article by C. Epstein and G. Henkin [EH] for the latest results and further references on preservation of embeddability for deformed CR structures on pseudoconvex hypersurfaces.

Version of Theorem 1 for CR submanifolds of codimension 1 (real hypersurfaces) is a corollary of a theorem proved by R.S. Hamilton [Ha2] provided that an additional assumption is satisfied: a deformation of the original CR structure is defined as a restriction of a deformation of a complex structure on some complex manifold $\mathbf{X}$ of which $\mathbf{M}$ is the boundary.

This property was proved by G.K. Kiremidjian [Ki] for a manifold $\mathbf{X}$, satisfying cohomological condition $H_c^2(\mathbf{X}, T'(\mathbf{X})) = 0$ and for deformations of a CR structure on a 2-pseudoconcave hypersurface $\mathbf{M}$, that preserve the underlying contact structure. Here, however, we do not use the extension property and work on $\mathbf{M}$ using constructions from [Ku], [Ki] and [Ha2].

In [Ha2] Hamilton applied estimates of Kohn and Nirenberg from [KN] and a special version of the Nash-Moser implicit function theorem for elliptic nonlinear complexes. Our approach to the proof of Theorem 1 basically follows S. Webster's scheme from [W2]. It consists of application of $C^k$-estimates for solutions of $\bar{\partial}_\mathbf{M}$-equation from [P4] and Zehnder's version (in Hamilton's terminology [Ha1]) of the Nash-Moser implicit function theorem. Necessary modifications are twofold. Firstly we use the language of differential forms (introduced by M. Kuranishi in [Ku]) instead of the language of vector fields as in [W2]. Secondly, the Webster's scheme is significantly simplified in our case because of the absence of "boundary terms" in the global homotopy formula.

The Nash-Moser iteration process is based on solvability with estimates of the linearized problem. In this paper we prove a global homotopy formula with *tame* $C^k$-estimates (cf. [Ha1]) for operator $\bar{\partial}_{\mathbf{M}}$ on a family of compact submanifolds $\mathbf{M}$, sufficiently close to a fixed submanifold $\mathbf{M}_0$. In the theorem below we formulate this solvability result.

**Theorem 2.** *Let $\mathbf{M}_0 \subset \mathbf{G}$ be a compact, generic, regular 3-pseudoconcave $C^p$ submanifold and let $\mathcal{L}$ be a holomorphic vector bundle on $\mathbf{G}$. Let the condition $\dim H^1\left(\mathbf{M}_0, \mathcal{L}|_{\mathbf{M}_0}\right) = 0$ be satisfied. Then for any $s, k \in \mathbb{Z}$ such that $k \geq 3$, $s \leq k$, $k + s \leq p - 5$ there exists $\delta > 0$ such that for any $C^p$ embedding*

$$\mathcal{E} : \mathbb{M} \to \mathbf{M} \subset \mathbf{G}$$

*with $|\mathcal{E}|_k < \delta$ there exist linear bounded operators*

$$\mathbf{P}_{\mathbf{M}} : C^{p-2}_{(0,1)}\left(\mathbf{M}, \mathcal{L}|_{\mathbf{M}}\right) \to C^{p-5}_{(0,0)}\left(\mathbf{M}, \mathcal{L}|_{\mathbf{M}}\right),$$

$$\mathbf{Q}_{\mathbf{M}} : C^{p-2}_{(0,2)}\left(\mathbf{M}, \mathcal{L}|_{\mathbf{M}}\right) \to C^{p-5}_{(0,1)}\left(\mathbf{M}, \mathcal{L}|_{\mathbf{M}}\right),$$

*such that for any differential form $h \in C^{p-2}_{(0,1)}\left(\mathbf{M}, \mathcal{L}|_{\mathbf{M}}\right)$ equality:*

$$h = \bar{\partial}_{\mathbf{M}} \mathbf{P}_{\mathbf{M}}(h) + \mathbf{Q}_{\mathbf{M}}(\bar{\partial}_{\mathbf{M}} h), \tag{2}$$

*and estimates*

$$|\mathbf{P}_{\mathbf{M}}(h)|_k, \ |\mathbf{Q}_{\mathbf{M}}(h)|_k \leq C(k)\left(1 + |\rho|_{k+5}\right)^{P(k)} |h|_{k+3}, \tag{3}$$

$$|\mathbf{P}_{\mathbf{M}}(h)|_{k+s}, \ |\mathbf{Q}_{\mathbf{M}}(h)|_{k+s} \leq C(k)\left(1 + |\rho|_{k+5}\right)^{P(k)} \left[\left(1 + |\rho|_{k+s+5}\right) |h|_{k+3} + |h|_{k+s+3}\right],$$

*hold, where $|\rho|_s = \max_{1 \leq i \leq N, 1 \leq l \leq m}\{|\rho_{i,l}|_s\}$, $\rho_{i,l}$ are defining functions of $\mathbf{M}$, and $P(k)$ is a polynomial in $k$.*

Analogous estimates were proved by S.Webster in [W1] for an integral formula on a special small neighborhood of a point on a strictly pseudoconvex hypersurface. Similar estimates were used by D.Catlin in [Ca] in his proof of extendability and local embeddability of CR structures. Our proof of existence of a tame homotopy formula basically consists of two parts: "local" and "global". We use local estimates from [P4] and complement them with consideration of "adjusted pairs" of neighborhoods, with motivation coming from the proof of local solvability of $\bar{\partial}_{\mathbf{M}}$-equation in [AiH]. The global part relies on a "soft analysis" of the deformation of $\bar{\partial}_{\mathbf{M}}$-complex and is based on the uniqueness theorem of S.Baouendi and F.Treves [BT] (cf. also [AiH]) for extensions of CR-functions from generic CR submanifolds.

Author thanks G. Henkin for helpful discussions.

## 2. Local solvability of $\bar{\partial}_{\mathbf{M}}$-equation with estimates.

We start with the discussion of local solvability of the $\bar{\partial}_{\mathbf{M}}$ - equation on pseudoconcave CR manifolds. Before constructing local solution formulas we introduce necessary notations and definitions.

For a vector-valued $C^1$-function $\eta = (\eta_1, \ldots, \eta_n)$ we will use the notation:

$$\omega'(\eta) = \sum_{k=1}^{n}(-1)^{k-1}\eta_k \wedge_{j \neq k} d\eta_j, \quad \omega(\eta) = \wedge_{j=1}^{n} d\eta_j.$$

If $\eta = \eta(\zeta, z, t)$ is a $C^1$-function of $\zeta \in \mathbb{C}^n, z \in \mathbb{C}^n$ and of a real parameter $t \in \mathbb{R}^l$, satisfying the condition

$$\sum_{k=1}^{n} \eta_k(\zeta, z, t) \cdot (\zeta_k - z_k) = 1 \tag{4}$$

then

$$d\omega'(\eta) \wedge \omega(\zeta) \wedge \omega(z) = 0$$

or, separating differentials,
$$d_t\omega'(\eta) + \bar{\partial}_\zeta\omega'(\eta) + \bar{\partial}_z\omega'(\eta) = 0. \tag{5}$$

Also, if $\eta(\zeta, z, t)$ satisfies (4) then the differential form $\omega'(\eta) \wedge \omega(\zeta) \wedge \omega(z)$ can be represented as:
$$\sum_{r=0}^{n-1} \omega'_r(\eta) \wedge \omega(\zeta) \wedge \omega(z), \tag{6}$$

where $\omega'_r(\eta)$ is a differential form of the order $r$ in $d\bar{z}$ and respectively of the order $n-r-1$ in $d\bar{\zeta}$ and $dt$. From (5) and (6) follow equalities:
$$d_t\omega'_r(\eta) + \bar{\partial}_\zeta\omega'_r(\eta) + \bar{\partial}_z\omega'_{r-1}(\eta) = 0 \qquad (r = 1, \ldots, n), \tag{7}$$

and
$$\omega'_r(\eta) = \frac{1}{(n-r-1)!r!} \mathrm{Det}\Big[\eta, \overbrace{\bar{\partial}_z\eta}^{r}, \overbrace{\bar{\partial}_{\zeta,t}\eta}^{n-r-1}\Big], \tag{8}$$

where the determinant is calculated by the usual rules but with external products of elements and the position of the element in the external product is defined by the number of its column.

Let $\mathcal{U}$ be an open neighborhood in $\mathbf{G}$ and $U = \mathcal{U} \cap \mathbf{M}$. We call a vector function
$$P(\zeta, z) = (P_1(\zeta, z), \ldots, P_n(\zeta, z)) \quad \text{for} \quad (\zeta, z) \in (\mathcal{U} \setminus U) \times U$$

by strong $\mathbf{M}$-barrier for $\mathcal{U}$ if there exists $C > 0$ such that the inequality:
$$|\Phi(\zeta, z)| > C \cdot \Big(\rho(\zeta) + |\zeta - z|^2\Big) \tag{9}$$

holds for $(\zeta, z) \in (\mathcal{U} \setminus U) \times U$, where
$$\Phi(\zeta, z) = \langle P(\zeta, z), \zeta - z \rangle = \sum_{i=1}^n P_i(\zeta, z) \cdot (\zeta_i - z_i).$$

If $\mathcal{U}$ is small enough, then as in (1) we may assume that $U = \mathcal{U} \cap \mathbf{M}$ is a set of common zeros of smooth functions $\{\rho_k, \ k = 1, \ldots, m\}$. Then, using the q-concavity of $\mathbf{M}$ and applying Kohn's lemma to the set of functions $\{\rho_k\}$ we can construct a new set of functions $\tilde{\rho}_1, \ldots, \tilde{\rho}_m$ of the form:
$$\tilde{\rho}_k(z) = \rho_k(z) + A \cdot \left(\sum_{i=1}^m \rho_i^2(z)\right),$$

with large enough constant $A > 0$ and such that for any $z \in \mathbf{M}$ there exist an open neighborhood $U \ni z$ and a family $E_{q+m}(\theta, z)$ of $q+m$ dimensional complex linear subspaces in $\mathbb{C}^n$ smoothly depending on $(\theta, z) \in \mathbb{S}^{m-1} \times U$ and such that $-\mathcal{L}_z\tilde{\rho}_\theta$ is strictly negative on $E_{q+m}(\theta, z)$ with all negative eigenvalues not exceeding some $c < 0$.

To simplify notations we will assume that the functions $\rho_1, \ldots, \rho_m$ already satisfy this condition.

Let $E^\perp_{n-q-m}(\theta, z)$ be the family of $n-q-m$ dimensional subspaces in $T(\mathbf{G})$ orthogonal to $E_{q+m}(\theta, z)$ and let
$$a_j(\theta, z) = (a_{j1}(\theta, z), \ldots, a_{jn}(\theta, z)) \text{ for } j = 1, \ldots, n-q-m$$

be a set of $C^2$-smooth vector functions representing an orthonormal basis in $E^\perp_{n-q-m}(\theta, z)$.

Defining for $(\theta, z, w) \in \mathbb{S}^{m-1} \times \mathcal{U} \times \mathbb{C}^n$
$$A_j(\theta, z, w) = \sum_{i=1}^n a_{ji}(\theta, z) \cdot w_i, \quad (j = 1, \ldots, n-q-m)$$

we construct the form

$$\mathcal{A}(\theta, z, w) = \sum_{j=1}^{n-q-m} A_j(\theta, z, w) \cdot \bar{A}_j(\theta, z, w)$$

such that the Hermitian form

$$\mathcal{L}_z \rho_\theta(w) + \mathcal{A}(\theta, z, w)$$

is strictly positive definite in $w$ for $(\theta, z) \in \mathbb{S}^{m-1} \times \mathcal{U}$.

Then we define for $\zeta, z \in (\mathcal{U} \setminus U) \times U$:

$$
\begin{aligned}
Q_i^{(k)}(\zeta, z) &= -\frac{\partial \rho_k}{\partial \zeta_i}(z) - \frac{1}{2} \sum_{j=1}^n \frac{\partial^2 \rho_k}{\partial \zeta_i \partial \zeta_j}(z)(\zeta_j - z_j), \\
F^{(k)}(\zeta, z) &= \langle Q^{(k)}(\zeta, z), \zeta - z \rangle, \\
P_i(\zeta, z) &= \sum_{k=1}^m \theta_k(\zeta) \cdot Q_i^{(k)}(\zeta, z) + \sum_{j=1}^{n-q-m} a_{ji}(\theta(\zeta), z) \cdot \bar{A}_j(\theta(\zeta), z, \zeta - z), \\
\Phi(\zeta, z) &= \langle P(\zeta, z), \zeta - z \rangle = \sum_{k=1}^m \theta_k(\zeta) \cdot F^{(k)}(\zeta, z) + \mathcal{A}(\theta(\zeta), z, \zeta - z)
\end{aligned}
\tag{10}
$$

with

$$\theta_k(\zeta) = -\frac{\rho_k(\zeta)}{\rho(\zeta)} \quad \text{for } k = 1, \ldots, m.$$

To prove that $P_i(\zeta, z)$ is a strong **M**-barrier for some $U \ni z$ we consider the Taylor expansion of $\rho_k$ for $k = 1, \ldots, m$:

$$\rho_k(\zeta) = \rho_k(z) - 2\mathrm{Re} F^{(k)}(\zeta, z) + \mathcal{L}_z \rho_k(\zeta - z) + O(|\zeta - z|^3).$$

Then we obtain for some $U$ and $(\zeta, z) \in (U \setminus (U \cap \mathbf{M})) \times (U \cap \mathbf{M})$:

$$
\begin{aligned}
\mathrm{Re}\Phi(\zeta, z) &= \sum_{k=1}^m \theta_k(\zeta) \cdot \mathrm{Re} F^{(k)}(\zeta, z) + \sum_{j=1}^{n-q-m} A_j(\theta(\zeta), \zeta, z) \cdot \bar{A}_j(\theta(\zeta), \zeta, z) \\
&= \rho(\zeta) + \mathcal{L}_z \rho_\theta(\zeta - z) + \mathcal{A}(\theta(\zeta), z, \zeta - z) + O(|\zeta - z|^3),
\end{aligned}
\tag{11}
$$

which implies the existence of an open neighborhood $U \ni z$ in $\mathbb{C}^n$, satisfying (9).

In what follows we will use notations

$$A_j(\zeta, z) := A_j(\theta(\zeta), z, \zeta - z) \text{ and } \mathcal{A}(\zeta, z) := \mathcal{A}(\theta(\zeta), z, \zeta - z).$$

In the lemma below we prove the existence of specially embedded convex neighborhoods in **G**.

**Lemma 2.1.** *Let* **M** *be a generic, regular $q$-pseudoconcave CR-submanifold of the class $C^p$ in $\mathcal{U}$ with $\mathcal{U}$ being a ball in a Riemannian metric on* **G** *centered at $z_0 \in$* **M**. *Then there exist a smaller ball $\mathcal{V} \subset \mathcal{U}$ centered at $z_0$ and $c > 0$ such that the barrier function $\Phi(\zeta, z)$, defined in (10), satisfies*

$$\inf_{(\mathbf{M} \cap b\mathcal{U}) \times \mathcal{V}} \mathrm{Re}\Phi(\zeta, z) > c. \tag{12}$$

**Proof.** Using (10) we obtain

$$
\begin{aligned}
\mathrm{Re}\Phi(\zeta, z) &= \sum_{k=1}^m \left(-\frac{\rho_k(\zeta)}{\rho(\zeta)}\right) \cdot \mathrm{Re} F^{(k)}(\zeta, z) + \mathcal{A}(\theta(\zeta), z, \zeta - z) \\
&= \rho(\zeta) - \sum_{k=1}^m \rho_k(z) \left(-\frac{\rho_k(\zeta)}{\rho(\zeta)}\right) + \mathcal{L}_z \rho_\theta(\zeta - z) + \mathcal{A}(\theta(\zeta), z, \zeta - z) + O(|\zeta - z|^3).
\end{aligned}
$$

Using then strict positivity of the form $\mathcal{L}_z\rho_\theta(w) + \mathcal{A}(\theta, z, w)$ we obtain that for $|\zeta - z|$ small enough there exists a constant $C > 0$ such that for $\zeta \in \mathbf{M} \cap b\mathcal{U}$

$$\mathrm{Re}\Phi(\zeta, z) \geq C|\zeta - z|^2 - \rho(z).$$

Therefore, for $z$, satisfying $\rho(z) < \frac{1}{2}C \inf_{\zeta \in \mathbf{M} \cap b\mathcal{U}} |\zeta - z|^2$ we have

$$\mathrm{Re}\Phi(\zeta, z) \geq \frac{1}{2}C \inf_{\zeta \in \mathbf{M} \cap b\mathcal{U}} |\zeta - z|^2.$$

Selecting a sufficiently small ball, centered at $z_0$ we obtain the statement of the lemma. □

For an open neighborhood $\mathcal{U} \subset \mathbf{G}$ we define the tubular neighborhood $\mathcal{U}(\epsilon)$ of $\mathbf{M}$ in $\mathcal{U}$ as

$$\mathcal{U}(\epsilon) = \{z \in \mathcal{U} : \rho(z) < \epsilon\},$$

where $\rho(z) = \left(\sum_{k=1}^m \rho_k^2(z)\right)^{\frac{1}{2}}$. We define $U(\epsilon)$ as

$$U(\epsilon) = \{z \in \mathcal{U} : \rho(z) = \epsilon\}.$$

We also introduce a local extension operator of functions and forms from $U = \mathcal{U} \cap \mathbf{M}$ to $\mathcal{U}$

$$E : g \to E(g).$$

Assuming that locally manifold $\mathbf{M}$ in $\mathcal{U}$ with coordinates $z_j = x_j + iy_j, j = 1, \ldots, n$ is defined as

$$\mathcal{U} \cap \mathbf{M} = \{z \in \mathcal{U} : \rho_j(z) \equiv x_j - h_j(y_1, \ldots, y_m, z_{m+1}, \ldots, z_n) = 0, j = 1, \ldots, m\}, \tag{13}$$

we define for a function $g(y_1, \ldots y_m, z_{m+1}, \ldots, z_n)$ on $U$

$$E(g)(z_1, \ldots, z_n) = g(y_1, \ldots, y_m, z_{m+1}, \ldots, z_n),$$

extending a function identically with respect to $x_1, \ldots, x_m$. For a differential form

$$g = \sum_{I,J,K} g_{I,J,K} dy_I \wedge dz_J \wedge d\bar{z}_K$$

with multiindices $I \subset (1, \ldots, m), J, K \subset (m+1, \ldots, n)$ we define extension operator by extending coefficients as in the formula above.

We start the construction of a local solution operator for $\bar{\partial}_\mathbf{M}$ with the construction of a Cauchy-Fantappie formula for special concave/convex neighborhoods. Let $\{\mathcal{U}_j\}_1^l$, $j \in J = (1, \ldots, l)$ be a collection of balls in some open neighborhood $\mathcal{W} \subset \mathbf{G}$, defined as

$$\mathcal{U}_j = \{\zeta \in \mathbf{G} : \tau_j(\zeta) < 0\},$$

where $\tau_j$ are convex $C^\infty$-functions. We denote also

$$\begin{aligned} \mathcal{U}_J &= \cap_{j=1}^l \mathcal{U}_j \neq \emptyset, \\ \mathcal{U}_J(\epsilon) &= \{\zeta \in \mathcal{U}_J : \rho(\zeta) < \epsilon\}. \end{aligned} \tag{14}$$

We assume the following nondegeneracy condition for functions $\{\rho_j\}_1^m, \{\tau_j\}_1^l$

$$\mathrm{rank}\{\mathrm{grad}\rho_1(\zeta), \ldots, \mathrm{grad}\rho_m(\zeta), \mathrm{grad}\tau_{j_1}(\zeta), \ldots, \mathrm{grad}\tau_{j_k}(\zeta)\} = m + k$$

for

$$\zeta \in \overline{\mathcal{U}_J(\epsilon)} \cap \{\tau_{j_1}(\zeta) = \cdots = \tau_{j_k}(\zeta) = 0\}.$$

Then the boundary of $\mathcal{U}(\epsilon)$ is stratified into $C^p$-smooth pieces

$$U^0(\epsilon) = \{\zeta \in \mathcal{U} : \rho(\zeta) = \epsilon, \ \tau_j(\zeta) < 0 \text{ for } j = 1, \ldots, l\},$$

$$U_I(\epsilon) = \{\zeta \in \mathcal{U} : \rho(\zeta) \leq \epsilon, \ \tau_i(\zeta) = 0 \text{ for } i \in I \subset J, \ \tau_i(\zeta) < 0 \text{ for } i \notin I\},$$

and
$$U_I^0(\epsilon) = \{\zeta \in \mathcal{U} : \rho(\zeta) = \epsilon, \ \tau_i(\zeta) = 0 \text{ for } i \in I \subset J, \ \tau_i(\zeta) < 0 \text{ for } i \notin I\}.$$

We introduce barriers for functions $\tau_j$:
$$\Phi^{(j)}(\zeta, z) = \sum P_i^{(j)}(\zeta)(\zeta_i - z_i) = \sum \frac{\partial \tau_j}{\partial \zeta_i}(\zeta)(\zeta_i - z_i),$$

such that $\mathrm{Re}\Phi^{(j)}(\zeta, z) > 0$ for $z \in \mathcal{U}_j$ and $\zeta \in b\mathcal{U}_j = \{\zeta : \tau_j(\zeta) = 0\}$.

In the proposition below we describe a general Cauchy-Fantappie formula for domains with piecewise concave/convex boundaries. We do not provide the proof of this proposition since it is completely analogous to the proof of the Cauchy-Fantappie formula in a piecewise strictly pseudoconvex case (cf. [P1] and [RS]) with the only difference that some of the barrier functions are concave and some are convex.

**Proposition 2.2.** *Let $\mathbf{M} \subset \mathbf{G}$ be a generic, regular q-pseudoconcave CR manifold of the class $C^p$ with $p \geq 3$ and $\mathcal{U}_J$ - a neighborhood as in (14) with analytic coordinates $z_1, \ldots, z_n$. Then for $\epsilon > 0$ small enough and an arbitrary differential form $g \in C^1_{(0,r)}(\mathcal{U}_J(\epsilon))$ the following equality holds for $r = 1, ..., n - m$ and $z \in U_J = \mathcal{U}_J \cap \mathbf{M}$*
$$g(z) = \bar{\partial} J_r(\epsilon)(g)(z) + J_{r+1}(\epsilon)(\bar{\partial} g)(z) + N_r(\epsilon)(g)(z), \tag{15}$$

*where*
$$J_r(\epsilon)(g)(z) = T_r(\epsilon)(g)(z) + R_r(\epsilon)(g)(z), \tag{16}$$

$$N_r(\epsilon)(g)(z) = S_r(\epsilon)(g)(z) + H_r(\epsilon)(g)(z), \tag{17}$$

$$T_r(\epsilon)(g)(z) = (-1)^r \frac{(n-1)!}{(2\pi i)^n} \int_{\mathcal{U}(\epsilon)} g(\zeta) \wedge \omega'_{r-1}\left(\frac{\bar{\zeta} - \bar{z}}{|\zeta - z|^2}\right) \wedge \omega(\zeta)$$
$$+ (-1)^r \frac{(n-1)!}{(2\pi i)^n} \sum_{1 \leq |I| \leq l} \int_{U_I(\epsilon) \times \Delta_I^s} g(\zeta) \wedge \omega'_{r-1}\left((1 - \sum_{i \in I} t_i)\frac{\bar{\zeta} - \bar{z}}{|\zeta - z|^2} + \sum_{i \in I} t_i \frac{P^{(i)}(\zeta, z)}{\Phi^{(i)}(\zeta, z)}\right) \wedge \omega(\zeta),$$

$$R_r(\epsilon)(g)(z) = (-1)^r \frac{(n-1)!}{(2\pi i)^n} \int_{U^0(\epsilon) \times [0,1]} g(\zeta) \wedge \omega'_{r-1}\left((1-t)\frac{\bar{\zeta} - \bar{z}}{|\zeta - z|^2} + t\frac{P(\zeta, z)}{\Phi(\zeta, z)}\right) \wedge \omega(\zeta)$$
$$+ (-1)^r \frac{(n-1)!}{(2\pi i)^n} \sum_{1 \leq |I| \leq l} \int_{U_I^0(\epsilon) \times \Delta_I^{s+1}} g(\zeta) \wedge \omega'_{r-1}\left((1 - t_0 - \sum_{i \in I} t_i)\frac{\bar{\zeta} - \bar{z}}{|\zeta - z|^2}\right.$$
$$\left.+ t_0 \frac{P(\zeta, z)}{\Phi(\zeta, z)} + \sum_{i \in I} t_i \frac{P^{(i)}(\zeta, z)}{\Phi^{(i)}(\zeta, z)}\right) \wedge \omega(\zeta),$$

$$S_r(\epsilon)(g)(z) = (-1)^r \frac{(n-1)!}{(2\pi i)^n} \sum_{1 \leq |I| \leq l} \int_{U_I(\epsilon) \times \Delta_I^{s-1}} g(\zeta) \wedge \omega'_r\left(\sum_{i \in I} t_i \frac{P^{(i)}(\zeta, z)}{\Phi^{(i)}(\zeta, z)}\right) \wedge \omega(\zeta),$$

$$H_r(\epsilon)(g)(z) = (-1)^r \frac{(n-1)!}{(2\pi i)^n} \int_{U^0(\epsilon)} g(\zeta) \wedge \omega'_r\left(\frac{P(\zeta, z)}{\Phi(\zeta, z)}\right) \wedge \omega(\zeta)$$
$$+ (-1)^r \frac{(n-1)!}{(2\pi i)^n} \sum_{1 \leq |I| \leq l} \int_{U_I^0(\epsilon) \times \Delta_I^s} g(\zeta) \wedge \omega'_r\left((1 - \sum t_i)\frac{P(\zeta, z)}{\Phi(\zeta, z)} + \sum_{i \in I} t_i \frac{P^{(i)}(\zeta, z)}{\Phi^{(i)}(\zeta, z)}\right) \wedge \omega(\zeta),$$

*and*
$$\Delta_I^{s+1} = \left\{(t_0, t_{i_1}, \ldots, t_{i_s}) \in \mathbb{R}^{s+1} : t_0 + \sum t_{i_k} \leq 1\right\},$$
$$\Delta_I^s = \left\{(t_{i_1}, \ldots, t_{i_s}) \in \mathbb{R}^s : \sum t_{i_k} \leq 1\right\},$$

$$\Delta_I^{s-1} = \left\{ (t_{i_1}, \ldots, t_{i_s}) \in \mathbb{R}^s : \sum t_{i_k} = 1 \right\},$$

for $I = (i_1, \ldots, i_s)$.

For $r = 0$ the corresponding formula becomes

$$g(z) = J_1(\epsilon)(\bar{\partial}g)(z) + N_0(\epsilon)(g)(z), \tag{18}$$

with $J_1(\epsilon)$ and $N_0(\epsilon)$ as above. $\square$

A local almost homotopy formula from the proposition below will be used in the construction of necessary solution operators.

**Proposition 2.3.** *Let $\mathbf{M} \subset \mathbf{G}$ be a generic, regular $q$-pseudoconcave CR submanifold of the class $C^p$ with $p \geq 3$. Let $\mathcal{U}_J$ be an open neighborhood as in (14) with analytic coordinates $z_1, \ldots, z_n$ and $V$ a relatively compact subset of $U_J = \mathbf{M} \cap \mathcal{U}_J$.*

*Then for $r = 1, \ldots, q-1$, $k \leq p$, and any differential form $g \in C^1_{(0,r)}(\mathbf{M})$ the equality*

$$g(z) = \bar{\partial}_{\mathbf{M}} R^r_{\mathbf{M}}(g)(z) + R^{r+1}_{\mathbf{M}}(\bar{\partial}_{\mathbf{M}} g)(z) + H^r_{\mathbf{M}}(g)(z), \tag{19}$$

*holds for $z \in V$, where*

$$R^r_{\mathbf{M}}(g)(z)$$
$$= (-1)^r \frac{(n-1)!}{(2\pi i)^n} \cdot pr_{\mathbf{M}} \circ \lim_{\epsilon \to 0} \left[ \int_{U^0(\epsilon) \times [0,1]} E(g)(\zeta) \wedge \omega'_{r-1} \left( (1-t)\frac{\bar{\zeta} - \bar{z}}{|\zeta - z|^2} + t\frac{P(\zeta, z)}{\Phi(\zeta, z)} \right) \wedge \omega(\zeta) \right.$$
$$+ \sum_{1 \leq |I| \leq |J|} \int_{U^0_I(\epsilon) \times \Delta_I^{s+1}} E(g)(\zeta) \wedge \omega'_{r-1} \left( (1 - t_0 - \sum_{i \in I} t_i)\frac{\bar{\zeta} - \bar{z}}{|\zeta - z|^2} \right.$$
$$\left. + t_0 \frac{P(\zeta, z)}{\Phi(\zeta, z)} + \sum_{i \in I} t_i \frac{P^{(i)}(\zeta, z)}{\Phi^{(i)}(\zeta, z)} \right) \wedge \omega(\zeta) \bigg],$$

$$H^r_{\mathbf{M}}(g)(z) = (-1)^r \frac{(n-1)!}{(2\pi i)^n} \cdot pr_{\mathbf{M}} \circ \lim_{\epsilon \to 0} \left[ \int_{U^0(\epsilon)} E(g)(\zeta) \wedge \omega'_r \left( \frac{P(\zeta, z)}{\Phi(\zeta, z)} \right) \wedge \omega(\zeta) \right.$$
$$+ \sum_{1 \leq |I| \leq |J|} \int_{U^0_I(\epsilon) \times \Delta_I^s} E(g)(\zeta) \wedge \omega'_r \left( (1 - \sum t_i)\frac{P(\zeta, z)}{\Phi(\zeta, z)} + \sum_{i \in I} t_i \frac{P^{(i)}(\zeta, z)}{\Phi^{(i)}(\zeta, z)} \right) \wedge \omega(\zeta) \bigg],$$

*$E(g)$ is an extension of $g$ to $\mathcal{U}_J$, and $pr_{\mathbf{M}}$ denotes the operator of projection to the space of tangential differential forms on $\mathbf{M}$.*

*For $r = 0$ the corresponding formula becomes*

$$g(z) = R^1_{\mathbf{M}}(\bar{\partial}_{\mathbf{M}} g)(z) + H^0_{\mathbf{M}}(g)(z), \tag{20}$$

*with $R^1_{\mathbf{M}}$ and $H^0_{\mathbf{M}}$ as above.*

**Proof.** We obtain formula (19) as a limit of the formula (15) and formula (20) as a limit of (18) when $\epsilon \to 0$. In order to make such a conclusion it suffices to prove that

$$\begin{aligned} \lim_{\epsilon \to 0} T_r(\epsilon)(g)(z) = 0, \\ \lim_{\epsilon \to 0} S_r(\epsilon)(g)(z) = 0, \end{aligned} \tag{21}$$

for $z \in V$.

To prove the first equality of (21) for the first integral in the definition of $T_r(\epsilon)$ we use the lemma below, which is a part of Proposition 4 from [P3].

**Lemma 2.4.** *Let $f \in L^s_{(0,r)}(U)$ with $s \in [1, \infty]$. Then for a relatively compact subset $V \subset U$*

$$\| T_r(\epsilon)(f) \|_{L^s_{(0,r)}(V)} = \mathcal{O}(\epsilon \cdot \log \epsilon) \cdot \| f \|_{L^s_{(0,r)}(U)} . \tag{22}$$

□

To prove that the sum of integrals over $U_I(\epsilon) \times \Delta^s_I$ in the definition of $T_r(\epsilon)$ also tends to zero as $\epsilon \to 0$ we notice that the kernels of those integrals are uniformly bounded for $z \in V$ and $\zeta \in U_I(\epsilon)$, $\dim U_I(\epsilon) = 2n - s$, and

$$\mathrm{mes}_{2n-s} U_I(\epsilon) = \mathcal{O}(\epsilon^m).$$

The same argument, applied to the sum of integrals in the definition of $S_r(\epsilon)$, proves the second equality in (21).

□

In the next proposition we prove a local solvability result for $\bar{\partial}_{\mathbf{M}}$-closed forms.

**Proposition 2.5.** *Let $\mathbf{M}$ be a generic, regular q-pseudoconcave submanifold of the class $C^p$ in $\mathbf{G}$, $J = (1, \ldots, l)$, $\{\mathcal{U}_j\}^l_1$ - a collection of balls in $\mathbf{G}$, centered respectively at the points $z_j \in \mathbf{M}$ and let $\{\mathcal{W}_j \subset \mathcal{V}_j\}^l_1$ and $\{\mathcal{V}_j \subset \mathcal{U}_j\}^l_1$ be the balls, satisfying (12) for any $j$. Let $\mathcal{U}_J = \cap^l_{j=1}\mathcal{U}_j$, $U_J = \mathbf{M} \cap \mathcal{U}_J$, $\mathcal{V}_J = \cap^l_{j=1}\mathcal{V}_j$, $V_J = \mathbf{M} \cap \mathcal{V}_J$, $\mathcal{W}_J = \cap^l_{j=1}\mathcal{W}_j$, and $W_J = \mathbf{M} \cap \mathcal{W}_J$.*

*If $W_J \neq \emptyset$, then for $1 \leq r < q$, $s \leq k$, and $k + s \leq p - 3$ there exists an operator*

$$P^r_{\mathbf{M}} : C^k_{(0,r)}(U_J) \to C^k_{(0,r-1)}(W_J)$$

*such that*

$$|P^r_{\mathbf{M}}(h)|_k \leq C(k)(1 + |\rho|_{k+3})^{P(k)} |h|_k, \tag{23}$$

$$|P^r_{\mathbf{M}}(h)|_{k+s} \leq C(k)(1 + |\rho|_{k+3})^{P(k)} [(1 + |\rho|_{k+s+3}) |h|_k + |h|_{k+s}],$$

*and equality*

$$h(z) = \bar{\partial}_{\mathbf{M}} P^r_{\mathbf{M}}(h)(z)$$

*holds for any $z \in W_J$ and any $h \in C^k_{(0,r)}(U_J)$ satisfying $\bar{\partial}_{\mathbf{M}} h = 0$.*

We list below several lemmas, that will be used in the proof of Proposition 2.5.
The following lemma is basically a copy of Lemma 3.8 from [P4].

**Lemma 2.6.** *Let $r < q$. Then*

$$\omega'_r \left( \frac{P(\zeta, z)}{\Phi(\zeta, z)} \right) = 0. \tag{24}$$

□

The lemma below is based on Lemma 4.1.1 and Theorem 4.2 from [AiH].

**Lemma 2.7.** *Let $\{\mathcal{U}_j\}^l_1$ and $\{\mathcal{V}_j \subset \mathcal{U}_j\}^l_1$ be the same as in Proposition 2.5. Then for any $\bar{\partial}_{\mathbf{M}}$-closed form $h \in C^k_{(0,r)}(U_J)$ with $3 \leq k$ and $0 \leq r < q$ the form*

$$\lim_{\epsilon \to 0} \sum_{\substack{I=(i_1,\ldots,i_s) \\ 1 \leq s \leq l}} \int_{U^0_I(\epsilon) \times \Delta^s_I} E(h)(\zeta) \wedge \omega'_r \left( (1 - \sum t_i) \frac{P(\zeta, z)}{\Phi(\zeta, z)} + \sum_{i \in I} t_i \frac{P^{(i)}(\zeta, z)}{\Phi^{(i)}(\zeta, z)} \right) \wedge \omega(\zeta)$$

*is defined and $\bar{\partial}$-closed on $\mathcal{V}_J = \cap^l_{j=1}\mathcal{V}_j$.*

**Proof.** We consider the differential form

$$\alpha(\zeta, z) = E(h)(\zeta) \wedge \omega'_r \left( (1 - \sum t_j) \frac{P(\zeta, z)}{\Phi(\zeta, z)} + \sum_{j \in J} t_j \frac{P^{(j)}(\zeta, z)}{\Phi^{(j)}(\zeta, z)} \right) \wedge \omega(\zeta)$$

for $z \in U_J$ and $(\zeta, t) \in \mathcal{C}$, where $\mathcal{C}$ is the chain

$$\mathcal{C} = \bigcup_{\substack{I \subset J \\ 1 \leq |I| \leq |J|}} \mathcal{C}_I = \bigcup_{\substack{I \subset J \\ 1 \leq |I| \leq |J|}} \left( U^0_I(\epsilon) \times \Delta^s_I \right).$$

Using formulas

$$b(\mathcal{C}_1 \times \mathcal{C}_2) = (\mathcal{C}_1 \times b\mathcal{C}_2) \cup (-1)^{\dim \mathcal{C}_2} (b\mathcal{C}_1 \times \mathcal{C}_2),$$

$$bU^0_I(\epsilon) = \cup_{k \notin I} U^0_{(I \cup k)}(\epsilon),$$

$$b\Delta^s_I = \left( \cup_{k=1}^s (-1)^k \Delta^{s-1}_{I(k)} \right) \cup \Delta^{s-1}_I,$$

where $I(k) = (i_1, \ldots, \widehat{i_k}, \ldots, i_s)$, we obtain the following formula for the boundary of each component $\mathcal{C}_I = U^0_I(\epsilon) \times \Delta^s_I$ of $\mathcal{C}$

$$b\mathcal{C}_I = (-1)^s \left[ \left( \cup_{k \notin I} U^0_{(I \cup k)}(\epsilon) \right) \times \Delta^s_I \right] \bigcup \left[ \cup_{k=1}^s (-1)^k \left( U^0_I(\epsilon) \times \Delta^{s-1}_{I(k)} \right) \cup \left( U^0_I(\epsilon) \times \Delta^{s-1}_I \right) \right].$$

Using then the formula above we obtain

$$b\mathcal{C} = \cup_{j \in J} \left( U^0_j(\epsilon) \times [1] \right) \bigcup (-1) \cup_{j \in J} \left( U^0_j(\epsilon) \times [0] \right) \bigcup_{I \subset J} \left( U^0_I(\epsilon) \times \Delta^{s-1}_I \right).$$

Applying $\bar\partial_z$ to the integral in the Lemma for $z \in \mathcal{V}_J$, using equality (7) and the Stokes' formula we obtain

$$\bar\partial_z \left[ \sum_{\substack{I=(i_1,\ldots,i_s) \\ 1 \leq s \leq l}} \int_{U^0_I(\epsilon) \times \Delta^s_I} E(h)(\zeta) \wedge \omega'_r \left( (1 - \sum t_i) \frac{P(\zeta, z)}{\Phi(\zeta, z)} + \sum_{i \in I} t_i \frac{P^{(i)}(\zeta, z)}{\Phi^{(i)}(\zeta, z)} \right) \wedge \omega(\zeta) \right] \quad (25)$$

$$= \int_{\mathcal{C}} E(h)(\zeta) \wedge \bar\partial_z \omega'_r \left( (1 - \sum t_i) \frac{P(\zeta, z)}{\Phi(\zeta, z)} + \sum_{i \in I} t_i \frac{P^{(i)}(\zeta, z)}{\Phi^{(i)}(\zeta, z)} \right) \wedge \omega(\zeta)$$

$$= \int_{\mathcal{C}} E(h)(\zeta) \wedge \bar\partial_{\zeta,t} \omega'_{r+1} \left( (1 - \sum t_i) \frac{P(\zeta, z)}{\Phi(\zeta, z)} + \sum_{i \in I} t_i \frac{P^{(i)}(\zeta, z)}{\Phi^{(i)}(\zeta, z)} \right) \wedge \omega(\zeta)$$

$$= \int_{b\mathcal{C}} E(h)(\zeta) \wedge \omega'_{r+1} \left( (1 - \sum t_i) \frac{P(\zeta, z)}{\Phi(\zeta, z)} + \sum_{i \in I} t_i \frac{P^{(i)}(\zeta, z)}{\Phi^{(i)}(\zeta, z)} \right) \wedge \omega(\zeta)$$

$$- \int_{\mathcal{C}} \bar\partial_\zeta E(h)(\zeta) \wedge \omega'_{r+1} \left( (1 - \sum t_i) \frac{P(\zeta, z)}{\Phi(\zeta, z)} + \sum_{i \in I} t_i \frac{P^{(i)}(\zeta, z)}{\Phi^{(i)}(\zeta, z)} \right) \wedge \omega(\zeta)$$

$$= \sum_j \int_{U^0_j(\epsilon) \times [1]} E(h)(\zeta) \wedge \omega'_{r+1} \left( \frac{P^{(j)}(\zeta, z)}{\Phi^{(j)}(\zeta, z)} \right) \wedge \omega(\zeta)$$

$$- \sum_j \int_{U^0_j(\epsilon) \times [0]} E(h)(\zeta) \wedge \omega'_{r+1} \left( \frac{P(\zeta, z)}{\Phi(\zeta, z)} \right) \wedge \omega(\zeta)$$

$$+ \sum_{\substack{I \subset J \\ 2 \leq |I| \leq |J|}} \int_{U^0_I(\epsilon) \times \Delta^{s-1}_I} E(h)(\zeta) \wedge \omega'_{r+1} \left( \sum_{i \in I} t_i \frac{P^{(i)}(\zeta, z)}{\Phi^{(i)}(\zeta, z)} \right) \wedge \omega(\zeta)$$

$$-\sum_{\substack{I\subset J\\1\leq|I|}}\int_{U_I^0(\epsilon)\times\Delta_I^s}\bar{\partial}_\zeta E(h)(\zeta)\wedge\omega'_{r+1}\left((1-\sum t_i)\frac{P(\zeta,z)}{\Phi(\zeta,z)}+\sum_{i\in I}t_i\frac{P^{(i)}(\zeta,z)}{\Phi^{(i)}(\zeta,z)}\right)\wedge\omega(\zeta).$$

Kernels in the first and third sums of integrals in the right hand side of (25) are zeros, since the barriers $\Phi^{(l)}$ are holomorphic in $z$. Therefore these sums are equal to zero.

For the fourth sum we use estimates

$$\operatorname{mes}_{2n-s-1}U_I^0(\epsilon)=\mathcal{O}(\epsilon^{m-1}),$$

$$\left|\omega'_{r+1}\left((1-\sum t_i)\frac{P(\zeta,z)}{\Phi(\zeta,z)}+\sum_{i\in I}t_i\frac{P^{(i)}(\zeta,z)}{\Phi^{(i)}(\zeta,z)}\right)\wedge\omega(\zeta)\right|_{U_I^0(\epsilon)}=\mathcal{O}(\epsilon^{-m+1}),$$

and

$$\left|\bar{\partial}_\zeta E(h)|_{U_I^0(\epsilon)}\right|_{U_I^0(\epsilon)}=\mathcal{O}(\epsilon)$$

for $h\in C^k_{(0,r)}(U_J)$ such that $k\geq 3$ and $\bar{\partial}_{\mathbf{M}}h|_U=0$, and obtain that this sum tends to zero as $\epsilon\to 0$.

For the second sum in the right hand side of (25) in the case $r<q-1$ we can apply Lemma 2.6 and obtain that this sum is zero. The case $r=q-1$ is more delicate, and we use lemma below, which is a specialization of Lemma 4.1.1 from [AiH] to domain $\mathcal{U}_J$.

**Lemma 2.8.** *Let $\mathbf{M}$, $\{\mathcal{U}_j\}_1^l$ and $\{\mathcal{V}_j\}_1^l$ be as in Proposition 2.5.*
*Then for any $\bar{\partial}_{\mathbf{M}}$-closed form $h\in C^k_{(0,q)}(U)$ with $k\geq 3$ and $z\in\mathcal{V}_J$ we have*

$$\lim_{\epsilon\to 0}\sum_j\int_{U_j^0(\epsilon)}E(h)(\zeta)\wedge\omega'_q\left(\frac{P(\zeta,z)}{\Phi(\zeta,z)}\right)\wedge\omega(\zeta)=0.$$

□

Summarizing all the estimates for the terms of the right hand side of (25) we obtain the statement of Lemma 2.7.

□

**Proof of Proposition 2.5.** We construct operators $P_{\mathbf{M}}^r$ taking formula (19) as a starting point. From this formula and Lemma 2.6 it follows that for $h\in C^k_{(0,r)}(U_J)$ satisfying $\bar{\partial}_{\mathbf{M}}h=0$ we have for any $z\in V_J$

$$h(z)=\bar{\partial}_{\mathbf{M}}R_{\mathbf{M}}^r(h)(z)+g(z),$$

with

$$g(z)=\lim_{\epsilon\to 0}\sum_{1\leq|I|\leq|J|}\int_{U_I^0(\epsilon)\times\Delta_I^s}E(h)(\zeta)\wedge\omega'_r\left((1-\sum t_i)\frac{P(\zeta,z)}{\Phi(\zeta,z)}+\sum_{i\in I}t_i\frac{P^{(i)}(\zeta,z)}{\Phi^{(i)}(\zeta,z)}\right)\wedge\omega(\zeta)\quad(26)$$

defined on $\mathcal{V}_J$ and satisfying equation $\bar{\partial}g=0$.

To obtain necessary estimates for operator $R_{\mathbf{M}}^r$ we use the following lemma, which is a combination of Propositions 3.1 and 3.8 from [P4].

**Lemma 2.9.** *Let $\mathbf{M}$ be a generic, regular $q$-pseudoconcave CR submanifold of the class $C^p$ in $\mathcal{U}$. Let $s,k\in\mathbb{Z}$ be such that $s\leq k$ and $k+s\leq p-3$.*
*Then $R_{\mathbf{M}}^r$ defined in (19) satisfies the following estimates*

$$|R_{\mathbf{M}}^r(h)|_k\leq C(k)(1+|\rho|_{k+3})^{P(k)}|h|_k,$$
$$|R_{\mathbf{M}}^r(h)|_{k+s}\leq C(k)(1+|\rho|_{k+3})^{P(k)}\left[|h|_{k+s}+(1+|\rho|_{k+s+3})|h|_k\right],\quad(27)$$

*with $P(k)$ a polynomial in $k$ and a constant $C(k)$ independent of $g$.*

□

To estimate the norm of $g$ we notice that the kernels of integrals in the definition of $g$ are nonsingular and uniformly bounded with derivatives of order $p - 2$ for $z \in \mathcal{V}_J$ and $\zeta \in U_J^0(\epsilon)$ for all $\epsilon > 0$. Straightforward differentiation of the formula (26) with respect to $z$ shows that

$$|g|_{\mathcal{V}_J, k+s} \leq C(k) \left(1 + |\rho|_{k+2}\right)^{P(k)} \left(1 + |\rho|_{k+s+2}\right) |h|_{U_J, 0}. \qquad (28)$$

Using then formula for solution of the equation

$$\bar{\partial} u = g$$

on the piecewise strictly pseudoconvex domain $\mathcal{V}_J$, analogous to formula (15) (cf. [P1],[RS]), we can represent solution $u$ as

$$u = u_{\mathcal{V}} + u_{b\mathcal{V}}, \qquad (29)$$

where $u_{\mathcal{V}}$ is a Martinelli-Bochner type integral over $\mathcal{V}_J$ and $u_{b\mathcal{V}}$ is a sum of integrals over the boundary $b\mathcal{V}_J$ with nonsingular kernels for $z \in \mathcal{W}$ and $\zeta \in b\mathcal{V}_J$. In this formula for the piecewise strictly pseudoconvex domain $\mathcal{V}_J$ part $N_r(\epsilon)$ will be absent, since the kernels will be holomorphic with respect to $z$ and part $J_{r+1}(\epsilon)(\bar{\partial} g)$ will be absent since $\bar{\partial} g = 0$.

We estimate two parts of (29) separately. For $u_{b\mathcal{V}}$ we have an estimate

$$|u_{b\mathcal{V}}|_{\mathcal{W}_J, k+s} \leq C(k) |g|_{\mathcal{V}_J, 0}.$$

For the Martinelli-Bochner integral using for example estimates from [Siu] we obtain

$$|u_{\mathcal{V}}|_{\mathcal{W}_J, k+s} \leq C(k) |g|_{\mathcal{V}_J, k+s}.$$

Defining then

$$P_{\mathbf{M}}^r(h) = R_{\mathbf{M}}^r(h) + u,$$

and combining the last two estimates with estimates (27) we obtain estimates (23). $\square$

In the proposition below we prove local extension with estimates for CR functions on a q-pseudoconcave CR submanifold.

**Proposition 2.10.** *Let $\mathbf{M}$ be a generic, regular q-pseudoconcave submanifold of the class $C^p$ in $\mathbf{G}$, $J = (1, \ldots, l)$, $\{\mathcal{U}_j\}_1^l$ - a collection of balls in $\mathbf{G}$, centered respectively at the points $z_j \in \mathbf{M}$ and let $\{\mathcal{V}_j \subset \mathcal{U}_j\}_1^l$ be the balls, satisfying (12) in each $\mathcal{U}_j$. Let $\mathcal{U}_J = \cap_{j=1}^l \mathcal{U}_j$, $U_J = \mathbf{M} \cap \mathcal{U}_J$, $\mathcal{V}_J = \cap_{j=1}^l \mathcal{V}_j$, and $V_J = \mathbf{M} \cap \mathcal{V}_J$.*

*Then there exists an extension operator*

$$P_{\mathbf{M}}^0 : \mathcal{O}_{\mathbf{M}}(U_J) \to \mathcal{O}(\mathcal{V}_J)$$

*acting on CR functions on $U_J$ and satisfying*

$$\left| P_{\mathbf{M}}^0(h) \right|_{k+s} \leq C(k) \left(1 + |\rho|_{k+2}\right)^{P(k)} \left(1 + |\rho|_{k+s+2}\right) |h|_0$$

*for any CR function $h \in \mathcal{O}_{\mathbf{M}}(U_J)$.*

**Proof.** We apply formula (20) to a CR function $h$ and obtain for $z \in V_J$

$$h(z) = H_{\mathbf{M}}^0(h)(z) = (-1)^r \frac{(n-1)!}{(2\pi i)^n} \cdot \mathrm{pr}_{\mathbf{M}} \circ \lim_{\epsilon \to 0} \left[ \int_{U^0(\epsilon)} E(h)(\zeta) \wedge \omega_0' \left( \frac{P(\zeta, z)}{\Phi(\zeta, z)} \right) \wedge \omega(\zeta) \right. \qquad (31)$$

$$\left. + \sum_{1 \leq |I| \leq |J|} \int_{U_I^0(\epsilon) \times \Delta_I^s} E(h)(\zeta) \wedge \omega_0' \left( (1 - \sum t_i) \frac{P(\zeta, z)}{\Phi(\zeta, z)} + \sum_{i \in I} t_i \frac{P^{(i)}(\zeta, z)}{\Phi^{(i)}(\zeta, z)} \right) \wedge \omega(\zeta) \right].$$

Using Lemma 2.6 we conclude that the first term of the formula above is equal to zero, and therefore, since kernels in the second term are nonsingular for $\zeta \in U_J^0(\epsilon)$ and $z \in \mathcal{V}_J$, function $H^0_{\mathbf{M}}(h)$ is defined on $\mathcal{V}_J$. Application of Lemma 2.7 shows that function $H^0_{\mathbf{M}}(h)$ is holomorphic. Straightforward differentiation of the formula (31) with respect to $z$ as in (28) shows that

$$\left|H^0_{\mathbf{M}}(h)\right|_{\mathcal{V}_J, k+s} \leq C(k) \left(1 + |\rho|_{k+2}\right)^{P(k)} \left(1 + |\rho|_{k+s+2}\right) |h|_{\mathcal{U}_J, 0}.$$

Defining then $P^0_{\mathbf{M}}(h) = H^0_{\mathbf{M}}(h)$ we obtain the statement of the Proposition. $\square$

## 3. From local to global estimates.

In this section we globalize estimates from the section 2. The following proposition (Proposition 4.5 from [P4]) provides a global formula on $\mathbf{M}_0$ that may be considered as an analogue of the Hodge-Kohn decomposition on $\mathbf{M}_0$.

**Proposition 3.1.** *Let $\mathbf{M}_0 \subset \mathbf{G}$ be a compact, generic, regular q-pseudoconcave submanifold of the class $C^p$ and let $\mathcal{L}$ be a holomorphic vector bundle on $\mathbf{G}$. Then for an arbitrary $r < q$ and $k \leq p - 4$ there exist a finite-dimensional operator*

$$K^r_{\mathbf{M}_0} : C^k_{(0,r)}\left(\mathbf{M}_0, \mathcal{L}|_{\mathbf{M}_0}\right) \to C^k_{(0,r)}\left(\mathbf{M}_0, \mathcal{L}|_{\mathbf{M}_0}\right),$$

*and operators $Q^r_{\mathbf{M}_0}$ and $Q^{r+1}_{\mathbf{M}_0}$ such that for $i = r, r+1$ the estimates*

$$\left|Q^i_{\mathbf{M}_0}(h)\right|_{k-1} \leq C(k) |h|_k, \tag{32}$$

*and equality*

$$h = \bar{\partial}_{\mathbf{M}_0} Q^r_{\mathbf{M}_0}(h) + Q^{r+1}_{\mathbf{M}_0}\left(\bar{\partial}_{\mathbf{M}_0} h\right) + K^r_{\mathbf{M}_0}(h) \tag{33}$$

*hold for any $h \in C^k_{(0,r)}\left(\mathbf{M}_0, \mathcal{L}|_{\mathbf{M}_0}\right)$.*

*If for $h \in C^k_{(0,r)}\left(\mathbf{M}_0, \mathcal{L}|_{\mathbf{M}_0}\right)$ there exists $g \in C^k_{(0,r-1)}\left(\mathbf{M}_0, \mathcal{L}|_{\mathbf{M}_0}\right)$ such that $\bar{\partial}_{\mathbf{M}_0} g = h$, then $K^r_{\mathbf{M}_0}(h) = 0$.* $\square$

We will reformulate the statement of Proposition 3.1 in terms of the Čech complex on $\mathbf{M}_0$. Let us consider the Čech complex for a bundle $\mathcal{L}$ and a covering $\{U\}$ of $\mathbf{M}_0$:

$$\Gamma\left(\mathbf{M}_0, \mathcal{O}_{\mathbf{M}_0}\right) \xrightarrow{\varrho} \mathcal{C}^0\left(\{U\}, \mathcal{O}_{\mathbf{M}_0}\right) \xrightarrow{\varrho} \mathcal{C}^1\left(\{U\}, \mathcal{O}_{\mathbf{M}_0}\right) \xrightarrow{\varrho} \mathcal{C}^2\left(\{U\}, \mathcal{O}_{\mathbf{M}_0}\right) \xrightarrow{\varrho} \cdots,$$

where $\Gamma\left(\mathbf{M}_0, \mathcal{O}_{\mathbf{M}_0}\right)$ is the space of CR sections of $\mathcal{L}$ over $\mathbf{M}_0$, $\mathcal{O}\left(U_I, \mathcal{L}|_{\mathbf{M}_0}\right)$ is the space of bounded CR sections of $\mathcal{L}|_{\mathbf{M}_0}$ over $U_I$,

$$\mathcal{C}^j\left(\{U\}, \mathcal{O}_{\mathbf{M}_0}\right) = \oplus_{|I|=j+1} \mathcal{O}\left(U_I, \mathcal{L}|_{\mathbf{M}_0}\right)$$

is the space of CR $j$-cochains in the covering $\{U_i = \mathbf{M}_0 \cap \mathcal{U}_i\}_1^N$ of $\mathbf{M}_0$ with $U_I = U_{i_0} \cap \cdots \cap U_{i_j}$ for $I = (i_0, \ldots, i_j)$ and the map $\varrho$ is defined as

$$\varrho(f)_{i_0 \ldots i_{j+1}} = \sum_{k=0}^{j+1} (-1)^k f_{i_0 \ldots \widehat{i_k} \ldots i_{j+1}}.$$

We call two coverings $\{\mathcal{V}_i \subset \mathcal{U}_i\}_1^N$ of some neighborhood of $\mathbf{M}_0$ adjusted if condition (12) is satisfied for every pair $\mathcal{V}_i \subset \mathcal{U}_i$.

In the next lemma we specialize the statement of Proposition 3.1 for 3-concave CR manifolds in terms of the Čech complex on $\mathbf{M}_0$.

**Lemma 3.2.** *Let $\mathbf{M}_0 \subset \mathbf{G}$ be a compact, generic, regular 3-pseudoconcave submanifold of the class $C^p$ with $p \geq 5$ and $\mathcal{L}$ be a holomorphic vector bundle on $\mathbf{G}$. Let $\{\mathcal{V}_i \subset \mathcal{U}_i\}_1^N$ be adjusted coverings of some neighborhood of $\mathbf{M}_0$. Let*

$$\mathcal{B}^2\left(\{U\}, \mathcal{O}_{\mathbf{M}_0}\right) = \left\{\alpha \in \mathcal{C}^2\left(\{U\}, \mathcal{O}_{\mathbf{M}_0}\right) : \exists \beta \in \mathcal{C}^1\left(\{U\}, \mathcal{O}_{\mathbf{M}_0}\right) \text{ such that } \varrho(\beta) = \alpha\right\}$$

*be the linear subspace of 2-coboundaries.*

*Then there exists an operator*

$$\mathcal{Q}_{\mathbf{M}_0}^2 : \mathcal{B}^2\left(\{U\}, \mathcal{O}_{\mathbf{M}_0}\right) \to \mathcal{C}^1\left(\{V\}, \mathcal{O}_{\mathbf{M}_0}\right)$$

*such that the estimate*

$$\left|\mathcal{Q}_{\mathbf{M}_0}^2(\alpha)\right|_{\{V\}} \leq C\, |\alpha|_{\{U\}}, \tag{34}$$

*and equality*

$$\alpha|_{\{V\}} = \varrho \circ \mathcal{Q}_{\mathbf{M}_0}^2(\alpha)$$

*hold for any $\alpha \in \mathcal{B}^2\left(\{U\}, \mathcal{O}_{\mathbf{M}_0}(\mathcal{L})\right)$.*

**Proof.** Let $\{\phi_i\}$ be a partition of unity, subordinate to the covering $\{V_i\}$ and therefore to $\{U_i\}$. For a covering $\{U_i\}$ we consider the double complex on $\mathbf{M}_0$

$$\begin{array}{ccccccccc}
\Gamma(\mathbf{M}_0, \mathcal{O}_{\mathbf{M}_0}) & \xrightarrow{\varrho} & \mathcal{C}^0(\{U\}, \mathcal{O}_{\mathbf{M}_0}) & \xrightarrow{\varrho} & \mathcal{C}^1(\{U\}, \mathcal{O}_{\mathbf{M}_0}) & \xrightarrow{\varrho} & \mathcal{C}^2(\{U\}, \mathcal{O}_{\mathbf{M}_0}) & \xrightarrow{\varrho} & \cdots \\
\downarrow & & \downarrow & & \downarrow & & \downarrow & & \\
\Gamma\left(\mathbf{M}_0, \Omega^{(0,0)}\right) & \xrightarrow{\varrho} & \mathcal{C}^0\left(\{U\}, \Omega^{(0,0)}\right) & \xrightarrow{\varrho} & \mathcal{C}^1\left(\{U\}, \Omega^{(0,0)}\right) & \xrightarrow{\varrho} & \mathcal{C}^2\left(\{U\}, \Omega^{(0,0)}\right) & \xrightarrow{\varrho} & \cdots \\
\downarrow \bar{\partial}_{\mathbf{M}_0} & & \downarrow \bar{\partial}_{\mathbf{M}_0} & & \downarrow \bar{\partial}_{\mathbf{M}_0} & & \downarrow \bar{\partial}_{\mathbf{M}_0} & & \\
\Gamma\left(\mathbf{M}_0, \Omega^{(0,1)}\right) & \xrightarrow{\varrho} & \mathcal{C}^0\left(\{U\}, \Omega^{(0,1)}\right) & \xrightarrow{\varrho} & \mathcal{C}^1\left(\{U\}, \Omega^{(0,1)}\right) & \xrightarrow{\varrho} & \mathcal{C}^2\left(\{U\}, \Omega^{(0,1)}\right) & \xrightarrow{\varrho} & \cdots \\
\downarrow \bar{\partial}_{\mathbf{M}_0} & & \downarrow \bar{\partial}_{\mathbf{M}_0} & & \downarrow \bar{\partial}_{\mathbf{M}_0} & & \downarrow \bar{\partial}_{\mathbf{M}_0} & & \\
\Gamma\left(\mathbf{M}_0, \Omega^{(0,2)}\right) & \xrightarrow{\varrho} & \mathcal{C}^0\left(\{U\}, \Omega^{(0,2)}\right) & \xrightarrow{\varrho} & \mathcal{C}^1\left(\{U\}, \Omega^{(0,2)}\right) & \xrightarrow{\varrho} & \mathcal{C}^2\left(\{U\}, \Omega^{(0,2)}\right) & \xrightarrow{\varrho} & \cdots \\
\downarrow \bar{\partial}_{\mathbf{M}_0} & & \downarrow \bar{\partial}_{\mathbf{M}_0} & & \downarrow \bar{\partial}_{\mathbf{M}_0} & & \downarrow \bar{\partial}_{\mathbf{M}_0} & & \\
\vdots & & \vdots & & \vdots & & \vdots & &
\end{array} \tag{35}$$

where

$$\Gamma\left(\mathbf{M}, \Omega^{(0,k)}\right) = C^{p-k}_{(0,k)}(\mathbf{M}, \mathcal{L}|_{\mathbf{M}})$$

is the space of differential forms of type $(0, k)$ on $\mathbf{M}$ with values in $\mathcal{L}|_{\mathbf{M}}$, and

$$\mathcal{C}^j\left(\{U\}, \Omega^{(0,k)}\right) = \oplus_{|I|=j+1} C^{p-k}_{(0,k)}(U_I, \mathcal{L}|_{\mathbf{M}})$$

is the space of $j$-cochains in the covering $\{U_i\}_1^N$.

We consider also operators

$$\chi : \mathcal{C}^k\left(\{U\}, \Omega^{(0,j)}\right) \to \mathcal{C}^{k-1}\left(\{U\}, \Omega^{(0,j)}\right),$$

defined by the formula

$$\chi(\alpha)_I = \sum_{l \notin I} \phi_l \alpha_{l \cup I},$$

and satisfying

$$\varrho \circ \chi(\alpha) = \alpha$$

for $\alpha$ such that $\varrho(\alpha) = 0$.

To an $\alpha \in \mathcal{B}^2(\{U\}, \mathcal{O}_{\mathbf{M}_0})$ we apply the *south-west diagram search* $\mathcal{SW}$, consisting of applications of vertical operators and operators $\chi$. Namely, we consider the 0-cochain of differential forms $\bar{\partial}_{\mathbf{M}_0}(\chi(\bar{\partial}_{\mathbf{M}_0}\chi(\alpha)))$, which is represented by a differential form $h = \mathcal{SW}(\alpha) \in C^{p-1}_{(0,2)}(\mathbf{M}_0, \mathcal{L}|_{\mathbf{M}_0})$ such that $\bar{\partial}_{\mathbf{M}_0} h = 0$ and $|h|_{p-1} \leq C |\alpha|_{\{U\}}$.

Since $\alpha$ satisfies $\alpha = \varrho(\beta)$ with some $\beta \in \mathcal{C}^1(\{U\}, \mathcal{O}_{\mathbf{M}_0})$, the south-west diagram search produces the form $h$ such that $h = \bar{\partial}_{\mathbf{M}_0} g$ with

$$g = \chi(\bar{\partial}_{\mathbf{M}_0}\chi(\alpha)) + \bar{\partial}_{\mathbf{M}_0}\chi(\beta - \chi(\alpha)) \in C^{p-1}_{(0,r-1)}(\mathbf{M}_0, \mathcal{L}|_{\mathbf{M}_0}).$$

Therefore, applying Proposition 3.1, we obtain representation

$$h = \bar{\partial}_{\mathbf{M}_0} Q^2_{\mathbf{M}_0}(h)$$

with the form $Q^2_{\mathbf{M}_0}(h)$, satisfying estimate

$$\left| Q^2_{\mathbf{M}_0}(h) \right|_{p-5} \leq C |\alpha|_{\{U\}}.$$

We consider then the 0-cochain of $\bar{\partial}_{\mathbf{M}_0}$-closed $(0,1)$ forms

$$\varrho\left(Q^2_{\mathbf{M}_0}(h)\right) - \chi(\bar{\partial}_{\mathbf{M}_0}\chi(\alpha))$$

and using Proposition 2.5 construct a 0-cochain of sections

$$\left\{ f_i = P^r_{\mathbf{M}_0}\left( \varrho\left(Q^2_{\mathbf{M}_0}(h)\right) - \chi(\bar{\partial}_{\mathbf{M}_0}\chi(\alpha)) \right) \right\} \in \mathcal{C}^0\left(\{V\}, \Omega^{(0,0)}\right),$$

satisfying estimate

$$|f_i|_{V_i, p-5} \leq C |\alpha|_{\{U\}},$$

and such that

$$\bar{\partial}_{\mathbf{M}_0} f_i = \left[ \varrho\left(Q^2_{\mathbf{M}_0}(h)\right) - \chi(\bar{\partial}_{\mathbf{M}_0}\chi(\alpha)) \right]\Big|_{V_i}.$$

Then 1-cochain $\mathcal{Q}^2(\alpha) = \varrho\{f_i\} + \chi(\alpha)$ satisfies equalities

$$\bar{\partial}_{\mathbf{M}_0} \mathcal{Q}^2(\alpha) = 0,$$

$$\varrho\left(\mathcal{Q}^2(\alpha)\right) = \alpha\big|_{\{V\}},$$

and estimate (34) and, therefore, defines the sought 1-cochain. $\square$

In the next proposition we prove the solvability of $\bar{\partial}_{\mathbf{M}}$ equation with tame estimates on manifolds $\mathbf{M}$ close to the fixed manifold $\mathbf{M}_0$.

**Proposition 3.3.** *Let $\mathbf{M}_0 \subset \mathbf{G}$ be a compact, generic, regular 3-pseudoconcave $C^p$ submanifold and $\{\mathcal{U}_i\}_1^N$ - a finite covering of some neighborhood of $\mathbf{M}_0$ in $\mathbf{G}$. Let*

$$B^k_{(0,r)}\left(\mathbf{M}, \mathcal{L}|_{\mathbf{M}}\right) = \left\{ h \in C^k_{(0,r)}\left(\mathbf{M}, \mathcal{L}|_{\mathbf{M}}\right) : \exists g \in C^k_{(0,r-1)}\left(\mathbf{M}, \mathcal{L}|_{\mathbf{M}}\right) \text{ such that } h = \bar{\partial}_{\mathbf{M}} g \right\}$$

*be a linear subspace in $C^k_{(0,r)}\left(\mathbf{M}, \mathcal{L}|_{\mathbf{M}}\right)$.*

*Then for $r = 1, 2$, $k \geq s > 3$, and $k + s \leq p - 4$ there exists $\delta > 0$ such that for any submanifold*

$$\mathcal{E} : \mathbb{M} \to \mathbf{M} \subset \mathbf{G}$$

*with defining functions $\{\rho_{i,l}\}_1^m$ in $\{\mathcal{U}_i\}_1^N$ and $|\mathcal{E}|_k < \delta$ there exist linear operators*

$$\mathbf{R}^r_{\mathbf{M}} : B^k_{(0,r)}\left(\mathbf{M}, \mathcal{L}|_{\mathbf{M}}\right) \to C^{k-1}_{(0,r-1)}\left(\mathbf{M}, \mathcal{L}|_{\mathbf{M}}\right)$$

*satisfying estimates*

$$|\mathbf{R}_{\mathbf{M}}^r(h)|_k \leq C(k)\left(1 + |\rho|_{k+4}\right)^{P(k)} |h|_{k+1}, \tag{36}$$

$$|\mathbf{R}_{\mathbf{M}}^r(h)|_{k+s} \leq C(k)\left(1 + |\rho|_{k+4}\right)^{P(k)} \left[\left(1 + |\rho|_{k+s+4}\right)|h|_{k+1} + |h|_{k+s+1}\right],$$

*and such that*

$$h = \bar{\partial}_{\mathbf{M}} \mathbf{R}_{\mathbf{M}}^r(h). \tag{37}$$

**Proof.** We describe the proof of the Proposition only for $r = 2$, since the proof for $r = 1$ is basically the same.

We consider adjusted coverings $\{\mathcal{U}_i^{(6)} \subset \cdots \subset \mathcal{U}_i^{(1)} = \mathcal{U}_i\}_1^N$ of some neighborhood $\mathcal{G}$ of $\mathbf{M}_0$ in $\mathbf{G}$. If $k > 3$, from the definition of adjusted coverings and from the constructions of Lemma 2.1 we conclude that there exists $\delta > 0$ such that selected coverings will be adjusted for manifolds $\mathbf{M}$ with $|\mathcal{E}|_k < \delta$. We consider double complexes analogous to (35) corresponding to those coverings with operators $\bar{\partial}_{\mathbf{M}}$ and $\varrho$, and fix a partition of unity subordinate to $\{\mathcal{U}^{(6)}\}$ and operators $\chi$ defined by this partition of unity.

To an arbitrary $h \in B_{(0,2)}^k\left(\mathbf{M}, \mathcal{L}|_{\mathbf{M}}\right)$ such that $h = \bar{\partial}_{\mathbf{M}} g$ we apply, using Propositions 2.5 and 2.10, the *north-east diagram search* $\mathcal{NE}$, consisting of the sequence of maps $\varrho$ and $P_{\mathbf{M}}^r$ for $r = 2, 1, 0$ and produce a Čech cochain $\mathcal{NE}(h) = \alpha \in \mathcal{C}^2\left(\{\mathcal{U}^{(4)}\}, \mathcal{L}\right)$ such that its restriction to $\mathbf{M}$ is a cocycle. To construct such a cochain we apply Proposition 2.5 to $h$ on $\mathcal{U}_i^{(1)}$ and consider forms $h_i^{(0,1)} = P_{\mathbf{M}}^2(h)$ on $U_i^{(2)}$ such that $\bar{\partial}_{\mathbf{M}} h_i^{(0,1)} = h$. Then we construct sections $h_{ij}^{(0,0)} = P_{\mathbf{M}}^1\left(h_i^{(0,1)} - h_j^{(0,1)}\right)$ on $U_{ij}^{(3)}$ such that $\bar{\partial}_{\mathbf{M}} h_{ij} = h_i - h_j$ and satisfying estimates

$$|h_{ij}|_k \leq C(k)\left(1 + |\rho|_{k+3}\right)^{P(k)} |h|_k, \tag{38}$$

$$|h_{ij}|_{k+s} \leq C(k)\left(1 + |\rho|_{k+3}\right)^{P(k)} \left[\left(1 + |\rho|_{k+s+3}\right)|h|_k + |h|_{k+s}\right]$$

for $k + s \leq p - 3$.

Finally, we define $\alpha = \{h_{ijk}\}$ on $\mathcal{U}_{ijk}^{(4)}$ with $h_{ijk} = P_{\mathbf{M}}^0(h_{ij} - h_{ik} + h_{jk})$, satisfying estimate

$$|\alpha|_{U^{(4)}} \leq C\left(1 + |\rho|_3\right)^P |h|_0.$$

Using uniqueness of holomorphic extension of CR functions (cf. [BT], [AiH]) we conclude that $\alpha$ is a cocycle on $\{\mathcal{U}^{(4)}\}$, and therefore, if $\mathbf{M} \subset \mathcal{G}$, the restriction of $\alpha$ to $\mathbf{M}_0$ is a cocycle as well.

Furthermore, since $h = \bar{\partial}_{\mathbf{M}} g$ with $g \in C_{(0,1)}^k\left(\mathbf{M}, \mathcal{L}|_{\mathbf{M}}\right)$, we apply the north-east diagram search to $g$ and construct a holomorphic cochain $\beta \in \mathcal{C}^1\left(\{\mathcal{U}^{(4)}\}, \mathcal{L}\right)$ such that $\alpha$ is a coboundary of $\beta$. Namely, consecutively using Proposition 2.5 we consider forms $h_i^{(0,1)} = P_{\mathbf{M}}^2(h)$ such that $\bar{\partial}_{\mathbf{M}} h_i^{(0,1)} = h$ on $U_i^{(2)}$ and then sections $f_i = P_{\mathbf{M}}^1(h_i - g)$ such that $\bar{\partial}_{\mathbf{M}} f_i = h_i - g$ on $U_i^{(3)}$. Then by construction we will have $\bar{\partial}_{\mathbf{M}}(f_i - f_j) = h_i^{(0,1)} - h_j^{(0,1)}$, and thus the holomorphic cochain $\{\beta_{ij} = P_{\mathbf{M}}^0(h_{ij} - f_i + f_j)\}$ on $\{\mathcal{U}_i^{(4)}\}$ will satisfy $\varrho(\beta|_{\mathbf{M}}) = \alpha|_{\mathbf{M}}$. Again using uniqueness of holomorphic extension for CR functions we obtain $\varrho(\beta) = \alpha$ on $\mathcal{G}$.

Considering then 2-cochain $\alpha|_{\mathbf{M}_0}$ on $\{\mathcal{U}_i^{(4)} \cap \mathbf{M}_0\}$ and applying Lemma 3.2, we obtain the 1-cochain $\mathcal{Q}^2(\alpha|_{\mathbf{M}_0})$ on $\{\mathcal{U}_i^{(5)} \cap \mathbf{M}_0\}$ satisfying estimate (34) and such that $\varrho\left(\mathcal{Q}^2(\alpha|_{\mathbf{M}_0})\right) = \alpha|_{\mathbf{M}_0}$. Another application of $P_{\mathbf{M}_0}^0$ produces the holomorphic 1-cochain $\gamma = P_{\mathbf{M}_0}^0\left(\mathcal{Q}^2(\alpha|_{\mathbf{M}_0})\right)$ on $\{\mathcal{U}^{(6)}\}$ such that $\varrho(\gamma)|_{\mathbf{M}_0} = \alpha|_{\mathbf{M}_0}$ and

$$|\gamma|_{\mathcal{U}^{(6)}} \leq C\left(1 + |\rho|_3\right)^P |h|_0. \tag{39}$$

Using uniqueness of holomorphic extension for CR functions we conclude that

$$\varrho(\gamma) = \alpha$$

on $\{\mathcal{U}^{(6)}\}$.

To construct operator $\mathbf{R}_{\mathbf{M}}^2$ we apply the south-west diagram search to cochain $\gamma$. Namely, we consider the 0-cochain of differential forms

$$\mathbf{R}_{\mathbf{M}}^2(h) = g_i^{(0,1)} = h_i^{(0,1)} - \bar{\partial}_{\mathbf{M}} \chi \left( h_{ij}^{(0,0)} - \gamma_{ij} \right) \text{ on } U_i^{(6)}.$$

Obviously condition $\bar{\partial}_{\mathbf{M}} g_i^{(0,1)} = h$ is satisfied. Also, by construction this cochain satisfies $\varrho\{g_i\} = 0$ and, therefore, defines a differential form on the whole $\mathbf{M}$. Estimates (36) follow from the construction of $\mathbf{R}_{\mathbf{M}}^2$ and estimates (38) and (39). □

**Proof of Theorem 2**.

To prove Theorem 2 we consider an arbitrary differential form $h \in C_{(0,1)}^{k+s+3} \left( \mathbf{M}, \mathcal{L}|_{\mathbf{M}} \right)$, and applying Proposition 3.3 to the differential form $\bar{\partial}_{\mathbf{M}} h \in C_{(0,2)}^{k+s+2} \left( \mathbf{M}, \mathcal{L}|_{\mathbf{M}} \right)$ construct a differential form $\mathbf{R}_{\mathbf{M}}^2(\bar{\partial}_{\mathbf{M}} h)$, satisfying

$$\bar{\partial}_{\mathbf{M}} \left( h - \mathbf{R}_{\mathbf{M}}^2(\bar{\partial}_{\mathbf{M}} h) \right) = 0$$

and

$$\left| \mathbf{R}_{\mathbf{M}}^2(\bar{\partial}_{\mathbf{M}} h) \right|_{k+1} \leq C(k) \left( 1 + |\rho|_{k+5} \right)^{P(k)} |h|_{k+3}, \tag{40}$$

$$\left| \mathbf{R}_{\mathbf{M}}^2(\bar{\partial}_{\mathbf{M}} h) \right|_{k+s+1} \leq C(k) \left( 1 + |\rho|_{k+5} \right)^{P(k)} \left[ \left( 1 + |\rho|_{k+s+5} \right) |h|_{k+3} + |h|_{k+s+3} \right].$$

From the assumption $\dim H^1 \left( \mathbf{M}_0, \mathcal{L}|_{\mathbf{M}_0} \right) = 0$ we obtain applicability of Proposition 3.3 to the differential form $h - \mathbf{R}_{\mathbf{M}}^2(\bar{\partial}_{\mathbf{M}} h) \in C_{(0,1)}^{k+s+1} \left( \mathbf{M}, \mathcal{L}|_{\mathbf{M}} \right)$. Applying Proposition 3.3 we construct

$$\mathbf{R}_{\mathbf{M}}^1 \left( h - \mathbf{R}_{\mathbf{M}}^2(\bar{\partial}_{\mathbf{M}} h) \right) \in C_{(0,0)}^{k+s} \left( \mathbf{M}, \mathcal{L}|_{\mathbf{M}} \right),$$

satisfying

$$\left| \mathbf{R}_{\mathbf{M}}^1 \left( h - \mathbf{R}_{\mathbf{M}}^2(\bar{\partial}_{\mathbf{M}} h) \right) \right|_k \leq C(k) \left( 1 + |\rho|_{k+5} \right)^{P(k)} |h|_{k+3}, \tag{41}$$

$$\left| \mathbf{R}_{\mathbf{M}}^1 \left( h - \mathbf{R}_{\mathbf{M}}^2(\bar{\partial}_{\mathbf{M}} h) \right) \right|_{k+s} \leq C(k) \left( 1 + |\rho|_{k+5} \right)^{P(k)} \left[ \left( 1 + |\rho|_{k+s+5} \right) |h|_{k+3} + |h|_{k+s+3} \right].$$

and obtain

$$h - \mathbf{R}_{\mathbf{M}}^2(\bar{\partial}_{\mathbf{M}} h) = \bar{\partial}_{\mathbf{M}} \mathbf{R}_{\mathbf{M}}^1 \left( h - \mathbf{R}_{\mathbf{M}}^2(\bar{\partial}_{\mathbf{M}} h) \right)$$

or

$$h = \bar{\partial}_{\mathbf{M}} \mathbf{P}_{\mathbf{M}}(h) + \mathbf{Q}_{\mathbf{M}}(\bar{\partial}_{\mathbf{M}} h),$$

with operators

$$\mathbf{P}_{\mathbf{M}}(h) = \mathbf{R}_{\mathbf{M}}^1 \left( h - \mathbf{R}_{\mathbf{M}}^2(\bar{\partial}_{\mathbf{M}} h) \right)$$

and

$$\mathbf{Q}_{\mathbf{M}}(h) = \mathbf{R}_{\mathbf{M}}^2(h)$$

satisfying estimates (3). □

## 4. Deformations of CR structures.

In this section we describe straightforward generalizations of some constructions from [Ku] and [Ki] to the case of higher codimension.

Let $\mathcal{E} : \mathbb{M} \to \mathbf{M}$ be a compact $C^p$-smooth CR submanifold of the form (1) in a complex $n$-dimensional manifold $\mathbf{G}$. For $k \leq p$ a $C^k$-subbundle $E'' \in \mathbf{C}T(\mathbb{M})$ is called an almost CR structure on $\mathbb{M}$ of finite distance to the fixed structure $T''(\mathbf{M})$ if it satisfies the following conditions:

$$(i) \quad E' \cap E'' = \{0\}, \quad \text{where} \quad E' = \overline{E''},$$

$$(ii) \quad \pi''|_{E''} : E'' \to T''(\mathbf{M}) \text{ is a } C^k \text{ isomorphism},$$

where $\pi''$ is the orthogonal projection of $\mathbf{C}T(\mathbb{M})$ onto $T''(\mathbf{M})$ along $T'(\mathbf{M}) \oplus \mathbf{N}$.

As follows from the definition of an almost CR structure of finite distance to the fixed CR structure it may be defined in terms of a differential form of type $(0,1)$ on $\mathbf{M}$ with values in the bundle $T'(\mathbf{M}) \oplus \mathbf{N}$. But for the purposes of comparison of the deformed CR structure on $\mathbb{M}$ with the complex structure of the ambient manifold $\mathbf{G}$, M. Kuranishi in [Ku] defines it in terms of a differential form of type $(0,1)$ on $\mathbf{M}$ with values in $T'(\mathbf{G})|_{\mathbf{M}}$.

Following [Ku] we consider the restriction to $T'(\mathbf{M}) \oplus \mathbf{N}$ of orthogonal projection $\pi' : \mathbf{C}T(\mathbf{G}) \to T'(\mathbf{G})$ along $T''(\mathbf{G})$. Since $T''(\mathbf{G}) \cap (T'(\mathbf{M}) \oplus \mathbf{N}) = \{0\}$ this restriction is an isomorphism. Therefore an isomorphism

$$\tau : T'(\mathbf{G}) \to T'(\mathbf{M}) \oplus \mathbf{N}$$

is well defined as the inverse to $\pi'|_{T'(\mathbf{M}) \oplus \mathbf{N}}$. The following proposition gives a description of almost CR structures of finite distance to the fixed induced CR structure in terms of differential forms.

**Proposition 4.1.** *If a $C^k$-subbundle $E''$ is an almost CR structure of finite distance to $T''(\mathbf{M})$, then there exists a unique differential form $\mu \in C^k_{(0,1)} \left( \mathbf{M}, T'(\mathbf{G})|_{\mathbf{M}} \right)$ such that*

$$E'' = \{Z - \tau \circ \mu(Z) : Z \in T''(\mathbf{M})\}. \tag{42}$$

*Conversely, if $\mu \in C^k_{(0,1)} \left( \mathbf{M}, T'(\mathbf{G})|_{\mathbf{M}} \right)$ and at each point $z \in \mathbf{M}$ the mapping $(\pi'' \circ (\overline{\tau \circ \mu})) \circ (\pi'' \circ (\overline{\tau \circ \mu})) : T''(\mathbf{M}) \to T''(\mathbf{M})$ does not have an eigenvalue 1, then equality (42) defines an almost CR structure of finite distance to $T''(\mathbf{M})$.*

We will denote by $T''_\mu(\mathbf{M})$ an almost CR structure defined on $\mathbb{M}$ by the form $\mu$ and equality (42).

An almost CR structure $E''$ is a CR structure if it is integrable, i.e. if for any two sections $L_1(z), L_2(z)$ of $E''$ the vector field $[L_1, L_2](z)$ is also a section of $E''$. In order to define integrability of $T''_\mu(\mathbf{M})$ in terms of $\mu$ we need some additional constructions.

Let $\mathcal{U} \in \mathbf{G}$ be a neighborhood of the point $z \in \mathbf{M}$. A set $\vartheta_1, \ldots, \vartheta_n \in C^k_1(\mathcal{U})$ of differential forms of degree 1 with $C^k$-coefficients is called a defining set for the almost CR structure $E''$ in $\mathcal{U}$ if for any $\zeta \in \mathbf{M} \cap \mathcal{U}$

(i) $\vartheta_i(X)(\zeta) = 0$ for $i = 1, \ldots, n$ and $X \in E''_\zeta$,

(ii) $\{\vartheta_i(\zeta)\}^n_1$ form a basis in the space of linear forms, satisfying (i).

We have the following description of formal integrability of the almost CR structure $E''$ in terms of its defining set.

**Proposition 4.2.** *An almost CR structure $E''$ of finite distance to $T''(\mathbf{M})$, is integrable if and only if for any $z \in \mathbf{M}$ there exist a neighborhood $\mathcal{U} \ni z$ in $\mathbf{G}$ and a defining set $\vartheta_1, \ldots, \vartheta_n \in C^k_1(\mathcal{U})$ for this almost CR structure in $\mathcal{U}$ such that*

$$d\vartheta_i = 0 \mod (\vartheta_1, \ldots, \vartheta_n) \text{ on } \mathbf{M} \cap \mathcal{U} \text{ for } (\alpha = 1, \ldots, n). \tag{43}$$

We consider special defining sets for an almost CR structure $T''_\mu(\mathbf{M})$.

**Proposition 4.3.** *Let $\{z_1, \ldots, z_n\}$ be a system of local analytic coordinates for $\mathbf{G}$ in a neighborhood $\mathcal{U} \in \mathbf{G}$ such that $\mathcal{U} \cap \mathbf{M} \neq \emptyset$ and let the form $\mu$ from Proposition 4.1 be represented in $\mathcal{U}$ as*

$$\mu = \sum_{i=1}^n \mu^i \frac{\partial}{\partial z_i},$$

*where $\mu^i = \sum_{j=1}^n \mu^i_j d\bar{z}_j$ are scalar $C^k$-differential forms of type $(0,1)$.*

*Then the set of differential forms*

$$\vartheta_i = dz_i + \mu^i \quad (i = 1, \ldots, n) \tag{44}$$

*is a defining set for $T''_\mu(\mathbf{M})$ in $\mathcal{U}$.*

Before applying Proposition 4.2 to the defining set (44) we will define special local generators in $T'(\mathbf{M})$ and $T''(\mathbf{M})$ and describe $\tau$ in terms of these generators.

Let us fix $z \in \mathbf{M}$ and neighborhoods: $\mathcal{U} \ni z$ in $\mathbf{G}$ with local analytic coordinates $\{z_1, \ldots, z_n\}$ and $U = \mathcal{U} \cap \mathbf{M}$. We consider differential forms on $\mathcal{U}$

$$\partial \rho_l = \sum_{i=1}^n \frac{\partial \rho_l}{\partial z_i} dz_i, \quad \bar{\partial} \rho_l = \sum_{i=1}^n \frac{\partial \rho_l}{\partial \bar{z}_i} d\bar{z}_i \quad \text{for} \quad l = 1, \ldots, m,$$

and vector fields on $U$:

$$P'_l = \sum_{i=1}^n p^i_l \frac{\partial}{\partial z_i} \in T'(\mathbf{G})|_U, \quad P''_l = \bar{P}'_l = \sum_{i=1}^n \bar{p}^i_l \frac{\partial}{\partial \bar{z}_i} \in T''(\mathbf{G})|_U \tag{45}$$

such that they satisfy conditions

$$(i) \quad \partial \rho_l(P'_s) = \bar{\partial} \rho_l(P''_s) = \begin{cases} 1 & \text{for } l = s, \\ 0 & \text{for } l \neq s, \end{cases} \tag{46}$$

$$(ii) \quad P'_l \perp T'(\mathbf{M}), \quad P''_l \perp T''(\mathbf{M}) \quad \text{for} \quad l = 1, \ldots, m.$$

Complex tangent vector fields on $U$

$$Z_i = \frac{\partial}{\partial z_i} - \sum_{l=1}^m \frac{\partial \rho_l}{\partial z_i} P'_l \quad \text{and} \quad \bar{Z}_i = \frac{\partial}{\partial \bar{z}_i} - \sum_{l=1}^m \frac{\partial \rho_l}{\partial \bar{z}_i} P''_l \quad \text{for} \quad i = 1, \ldots, n,$$

span $T'(\mathbf{M})$ and $T''(\mathbf{M})$ respectively, and satisfy the following conditions

$$\sum_{i=1}^n p^i_l Z_i = \sum_{i=1}^n \bar{p}^i_l \bar{Z}_i = 0 \quad \text{for} \quad l = 1, \ldots, m.$$

To describe the isomorphism $\tau : T'(\mathbf{G}) \to T'(\mathbf{M}) \oplus \mathbf{N}$ we notice at first that it follows from (46) that

$$d\rho_k(P'_l - P''_l) = 0 \quad \text{for} \quad k, l = 1, \ldots, m,$$

and, therefore, $P'_l - P''_l \in \mathbf{N} \subset \mathbf{C}T(\mathbf{M})$.

Using then that $\tau(P'_l)$ is uniquely defined by the conditions

$$\tau(P'_l) - P'_l \in T''(\mathbf{G}),$$

$$\tau(P'_l) \in \mathbf{C}T(\mathbf{M}),$$

we conclude that

$$\tau(P'_l) = P'_l - P''_l.$$

For vector fields $Z_i \in T'(\mathbf{M})$ we have

$$\tau(Z_i) = Z_i$$

and, therefore,

$$\frac{\partial^\tau}{\partial z_i} := \tau\left(\frac{\partial}{\partial z_i}\right) = \tau\left(Z_i + \sum_{l=1}^m \frac{\partial \rho_l}{\partial z_i} P'_l\right) \qquad (47)$$

$$= Z_i + \sum_{l=1}^m \frac{\partial \rho_l}{\partial z_i}(P'_l - P''_l) = \frac{\partial}{\partial z_i} - \sum_{l=1}^m \frac{\partial \rho_l}{\partial z_i} P''_l.$$

For a differentiable function $\tilde{h}$ on $\mathcal{U}$ we can express its differential in terms of $\bar{Z}^i, dz_i$ and $d\rho_l$ as

$$d\tilde{h} = \sum_i \frac{\partial \tilde{h}}{\partial z_i} dz_i + \sum_i \frac{\partial \tilde{h}}{\partial \bar{z}_i} d\bar{z}_i = \sum_i \frac{\partial \tilde{h}}{\partial \bar{z}_i} \bar{Z}^i + \sum_i \frac{\partial^\tau \tilde{h}}{\partial z_i} dz_i + \sum_{l=1}^m P''_l\left(\tilde{h}\right) d\rho_l. \qquad (48)$$

Restricting (48) onto $\mathbf{M}$ we obtain for a differentiable function $h = \tilde{h}\big|_{\mathcal{U} \cap \mathbf{M}}$

$$dh = \sum_i \frac{\partial \tilde{h}}{\partial \bar{z}_i} \bar{Z}^i + \sum_i \frac{\partial^\tau \tilde{h}}{\partial z_i} dz_i, \qquad (49)$$

where $\sum_i \bar{Z}^i (\partial/\partial \bar{z}_i) = \bar{\partial}_\mathbf{M}$ is the Cauchy-Riemann operator on $\mathbf{M}$.

The following proposition formulates integrability of an almost CR structure $T''_\mu(\mathbf{M})$ in terms of $\mu$.

**Proposition 4.4.** ([Ku], cf. also [Ki]). *An almost CR structure $T''_\mu(\mathbf{M})$ is integrable if and only if*

$$\Phi(\mu) = 0, \qquad (50)$$

*where $\Phi(\mu)$ is a $T'(\mathbf{G})$ - valued $(0,2)$ differential form on $\mathbf{M}$, locally defined as $\Phi(\mu) := \sum_i \Phi^i(\mu) \frac{\partial}{\partial z_i}$ with*

$$\Phi^i(\mu) := \bar{\partial}_\mathbf{M} \mu^i - \sum_{k,j} \left(\frac{\partial^\tau \mu^i_j}{\partial z_k}\right) \mu^k \wedge \bar{Z}^j$$

$$- \sum_j \mu^i_j \cdot \sum_{l=1}^m \left(\bar{\partial}_\mathbf{M} \bar{p}^j_l - \sum_k \frac{\partial^\tau \bar{p}^j_l}{\partial z_k} \mu^k\right) \wedge \mu(\rho_l)$$

*and*

$$\mu(\rho_l) := \sum_j \frac{\partial \rho_l}{\partial z_j} \mu^j.$$

**Proof.** We fix a neighborhood $\mathcal{U}$ in $\mathbf{G}$ with coordinates $z_i$ $(i = 1, \ldots, n)$ and denote $U = \mathcal{U} \cap \mathbf{M}$. Then differential form $\mu$ on $\mathcal{U}$ may be represented as $\tilde{\mu} = \sum_{i,j} \mu^i_j \frac{\partial}{\partial z_i} \otimes d\bar{z}_j$. Without loss of generality we may assume that the coefficients of this form satisfy the condition

$$\sum_j \mu^i_j \bar{p}^j_k = 0 \text{ for } z \in U, \text{ any } i, \text{ and } k = 1, \ldots, m. \qquad (51)$$

Really, if the differential form $\tilde{\mu}^i$ doesn't satisfy (51) we can modify $\tilde{\mu}^i$ by setting it equal to

$$\tilde{\mu}^{i,\prime} = \tilde{\mu}^i - \sum_{l=1}^m \left(\sum_j \mu^i_j \bar{p}^j_l\right) \bar{\partial} \rho_l.$$

Then the form $\tilde{\mu}^{i,\prime}$ will coincide with $\tilde{\mu}^i$ on $U$ and satisfy (51) because of (46).

To prove Proposition 4.4 it suffices, according to propositions 4.2 and 4.3, to prove that condition (50) is necessary and sufficient for the defining set (44) to satisfy:

$$d\vartheta_i = 0 \mod (\vartheta_1, \ldots, \vartheta_n), \quad (i = 1, \ldots, n),$$

or, equivalently,

$$d\mu^i = 0 \mod (\vartheta_1, \ldots, \vartheta_n), \quad (i = 1, \ldots, n).$$

Using equalities

$$\bar{Z}^j\Big|_U = d\bar{z}_j - \sum_{l=1}^m \bar{p}_l^j \bar{\partial}\rho_l\Big|_U = d\bar{z}_j + \sum_{l=1}^m \bar{p}_l^j \partial\rho_l\Big|_U \quad \text{for} \quad j = 1, \ldots, n,$$

condition (51), and applying (49) to $\mu_j^i$ and $\bar{p}_l^j$, we obtain:

$$d\mu^i = \sum_j d\mu_j^i \wedge \bar{Z}^j + \sum_j \mu_j^i d\bar{Z}^j \tag{52}$$

$$= \sum_{k,j} \frac{\partial \mu_j^i}{\partial \bar{z}_k} \bar{Z}^k \wedge \bar{Z}^j + \sum_{k,j} \frac{\partial^\tau \mu_j^i}{\partial z_k} dz_k \wedge \bar{Z}^j$$

$$+ \sum_j \mu_j^i \cdot \sum_{l=1}^m \left( \bar{\partial}_{\mathbf{M}} \bar{p}_l^j + \sum_k \frac{\partial^\tau \bar{p}_l^j}{\partial z_k} dz_k \right) \wedge \partial\rho_l.$$

Using then in (52) equalities

$$dz_k = -\mu^k \mod (\vartheta_1, \ldots, \vartheta_n) \quad \text{and} \quad \partial\rho_l = -\mu(\rho_l) \mod (\vartheta_1, \ldots, \vartheta_n),$$

we obtain that locally condition (50) is necessary and sufficient for integrability of $T''_\mu(\mathbf{M})$.

To prove that $\Phi(\mu)$ is a differential form on $\mathbf{M}$ we notice first that $d\mu^i$ is a differential form on $\mathbf{M}$ independent of the choice of $\rho_l$ and $p_l$. Secondly, the forms $dz_k$, $\partial\rho_l$ change by the same rules as respectively $\mu^k$, $\mu(\rho_l)$ under the change of local coordinates. $\square$

**Remark.** Following [Ku] (cf. also [Ki]) we can define operator $\bar{\partial}_{\mathbf{M}}^\mu$ which represents the $\bar{\partial}_{\mathbf{M}}$ operator in an almost CR structure, defined on $\mathbb{M}$ by the form $\mu$ using formulas

$$\bar{\partial}_{\mathbf{M}}^\mu f = \bar{\partial}_{\mathbf{M}} f - \sum_j \frac{\partial^\tau f}{\partial z_j} \mu^j \quad \text{for} \quad f \in C^p(\mathbf{M})$$

and

$$\bar{\partial}_{\mathbf{M}}^\mu \bar{Z}^j = \sum_k \mu^k \wedge \sum_{l=1}^m \frac{\partial \rho_l}{\partial z_k} \bar{\partial}_{\mathbf{M}}^\mu \bar{p}_k^j.$$

As it is pointed out in [Ku] a straightforward calculation leads to equality

$$\bar{\partial}_{\mathbf{M}}^\mu \circ \bar{\partial}_{\mathbf{M}}^\mu f = -\sum_i \frac{\partial^\tau f}{\partial z_i} \Phi(\mu)^i$$

showing exactness of the $\bar{\partial}_{\mathbf{M}}^\mu$ - complex for $\mu$ satisfying integrability condition (50).

## 5. Properties of embeddings.

In this section we describe properties of an embedding $\mathcal{F} : \mathbb{M} \to \mathbf{G}$ of a compact manifold $\mathbb{M}$ close to a fixed embedding $\mathcal{E}$ with $\mathbf{M}_0 = \mathcal{E}(\mathbb{M})$ being a $C^p$-smooth, generic, regular 3-pseudoconcave CR manifold. We describe this embedding in terms of a map $F = \mathcal{F} \circ \mathcal{E}^{-1} : \mathbf{M} \to \mathbf{G}$ close to identity and summarize the estimates needed in the proof of Theorem 1.

Without loss of generality we may assume that in every neighborhood $\mathcal{U}$ of a finite covering of $\mathbf{G}$ the map $F$ is a restriction of a map $F : \mathcal{U} \to \mathcal{V}$ defined by the functions

$$F_i(z) = z_i + f_i(z), \quad i = 1, \ldots, n, \tag{53}$$

with

$$|f|_k \equiv \sup_{z \in \mathcal{U}, |J| \leq k} \left| \partial^J f(z) \right| < \epsilon,$$

where $J = (j_1, \ldots, j_{2n})$ is a multiindex, $|J| = j_1 + \cdots + j_{2n}$, and $\partial^J \equiv \partial_{x_1}^{j_1} \cdots \partial_{x_{2n}}^{j_{2n}}$ with coordinates $\{x_j\}_1^{2n}$ such that $z_j = x_j + \sqrt{-1} x_{n+j}$.

For $\epsilon$ small enough according to the inverse function theorem there exists an inverse map $G : \mathcal{U}^* \to \mathcal{U}$ defined on a possibly smaller neighborhood $\mathcal{U}^*$ by the functions

$$G_i(z) = z_i + g_i(z), \quad i = 1, \ldots, n, \tag{54}$$

such that $g_i \in C^p(\mathcal{U}^*)$ and

$$F \circ G(z) = z. \tag{55}$$

The proposition below, which is a copy of Proposition 4.6 from [P4], provides necessary estimates for functions $g_i$. We use the following notations. We denote by $\mathbb{B}(r)$ the ball in $\mathbb{C}^n$ of radius $r$ centered at the origin and for a function $f : \mathbb{B}(r) \to \mathbb{C}$ and $k \leq p$ we denote

$$|f|_{r,k} \equiv \sup_{z \in \mathbb{B}(r), |J| \leq k} \left| \partial^J f(z) \right|.$$

For a vector function $f : \mathbb{B}(r) \to \mathbb{C}^n$ we denote

$$|f|_{r,k} \equiv \sup_{1 \leq i \leq n} |f_i|_{r,k}.$$

**Proposition 5.1.** ( [P4], cf.[W2] ) *Let $F_i(z) = z_i + f_i(z)$ for $i = 1, \ldots, n$, and let the functions $\{f_i\}_1^n \in C^p(\mathbb{B}(1))$ satisfy estimate $|f_i|_{1,1} < \epsilon$.*

*Then for small enough $\epsilon$ and fixed $s, k \in \mathbb{Z}$, such that $0 \leq s \leq k$, $k + s \leq p$ there exist a constant $C(k)$ and a set of functions*

$$\{g_i\}_1^n \in C^p(\mathbb{B}(1 - 2\epsilon))$$

*such that $G(z) \equiv z + g(z) \in \mathbb{B}(1)$ for $z \in \mathbb{B}(1 - 2\epsilon)$,*

$$F \circ G(z) = z$$

*is satisfied on $\mathbb{B}(1 - 2\epsilon)$, and*

$$|g_i|_{1-2\epsilon, k} \leq C(k) \cdot (1 + |f|_{1,k})^{P(k)} |f|_{1,k}, \tag{56}$$

$$|g_i|_{1-2\epsilon, k+s} \leq C(k) \cdot (1 + |f|_{1,k})^{P(k)} |f|_{1,k+s}$$

*with polynomial $P(k)$.*

□

Let now map $F$ be as in (53). Then according to Proposition 5.1 the part of the image of the manifold $\mathbf{M} = \{z : \rho_l(z) = 0, l = 1, \ldots, m\}$ under $F$ that lies in $\mathcal{U}^*$ is the manifold

$$\mathbf{M}^* = F(\mathbf{M}) = \{z : \rho_l(G(z)) = 0, l = 1, \ldots, m\},$$

defined by the functions $\rho_l^*(z) \equiv \rho_l(z + g(z))$. In two lemmas below we prove necessary estimates for vector fields $\{P_l'\}_1^m$ and $\{P_l'^{,*}\}_1^m$ satisfying conditions (46) on $\mathbf{M}$ and $\mathbf{M}^*$ respectively.

**Lemma 5.2.** *Let functions $\{\rho_l\}_1^m \in C^p(\mathbb{B}(1))$ have the form*

$$\rho_l(z) = y_l + r_l(z) \quad \text{for} \quad l = 1, \ldots, m \tag{57}$$

*in coordinates $z_1, \ldots, z_n$ with*

$$\partial r_l / \partial z_j(0) = \partial r_l / \partial \bar{z}_j(0) = 0 \text{ for } l = 1, \ldots, m \text{ and } j = 1, \ldots, n.$$

*Then for a small enough neighborhood $\mathcal{U} \ni 0$ there exist vector fields $\{P_l'\}_1^m$ satisfying (46) and estimates*

$$\left| P' \right|_k \leq C(k) \left( 1 + |\rho|_{k+1} \right)^{P(k)}, \tag{58}$$

$$|P'|_{k+s} \leq C(k)\,(1+|\rho|_{k+1})^{P(k)}\,(1+|\rho|_{k+s+1})$$

for $s, k \in \mathbb{Z}$, $s \leq k$, $s + k \leq p - 1$ with some polynomial $P(k)$.

**Proof.** If $\mathbf{M}$ is defined in $\mathcal{U}$ by the functions $\rho_l$, $l = 1, \ldots, m$ as in (57) then for a small enough neighborhood $\mathcal{U}$, denoting by $^\tau$ the transposition operator, we will have the following representation

$$\left[\frac{\partial \rho}{\partial z}\right]\left[\frac{\partial \rho}{\partial \bar{z}}\right]^\tau = [I + A] \tag{59}$$

with the matrix $A$ admitting estimates

$$|A|_0 < \delta \ll 1, \tag{60}$$

$$|A|_k \leq C(k)\,(1+|\rho|_{k+1})^{P(k)},$$

$$|A|_{k+s} \leq C(k)\,(1+|\rho|_{k+1})^{P(k)}\,(1+|\rho|_{k+s+1}).$$

Assuming then that the metric in $\mathcal{U}$ is the standard Hermitian metric in coordinates $z$, and choosing

$$[P'] = \left[\frac{\partial \rho}{\partial \bar{z}}\right]^\tau \cdot [I + A]^{-1}, \tag{61}$$

we obtain the vector fields, satisfying (46) and estimates (58). $\square$

**Lemma 5.3.** *Let the functions $\{\rho_l\}_1^m \in C^p(\mathbb{B}(1))$ be chosen as in (57) and let the functions $\{g_j\}_1^n \in C^p(\mathbb{B}(1))$ be such that $|g|_k < \epsilon < 1$ for some fixed $k$ such that $2 \leq k \leq p - 2$. Then for small enough $\epsilon$ and a neighborhood $\mathcal{U} \ni 0$ there exist vector fields*

$$P'^{,*}_l = \sum_{i=1}^n p_l^{*,i} \frac{\partial}{\partial z_i} \in C^{p-1}\left(T'(\mathbb{B}(0, 1 - C \cdot \epsilon))\right),$$

*satisfying (46) for $\{\rho_l^*\}_1^m$ and such that for $s, k \in \mathbb{Z}$, $s \leq k$, $s + k \leq p - 2$*

$$|P'^{,*}_l - P'_l|_k \leq C(k)\,(1+|\rho|_{k+1})^{P(k)}\,|g|_k, \tag{62}$$

$$|P'^{,*}_l - P'_l|_{k+s} \leq C(k)\,(1+|\rho|_{k+1})^{P(k)}\,[(1+|\rho|_{k+s+1})\,|g|_k + |g|_{k+s}],$$

*where $P(k)$ is a polynomial in $k$.*

**Proof.** Application of the chain rule to functions $\{\rho_l^*\}_1^m$ gives an estimate for $l = 1, \ldots, m$:

$$|\rho_l^* - \rho_l|_k = |\rho_l(z + g(z)) - \rho_l(z)|_k$$

$$= \left|\int_0^1 \langle \nabla \rho_l(z + tg(z)), g(z)\rangle dt\right|_k \leq C(k)\,(1+|\rho|_{k+1})^{P(k)}\,|g|_k \tag{63}$$

with a polynomial $P(k)$ and a constant $C(k)$ depending only on $k$.

Analogously, for $k + s \leq p - 2$ and $s \leq k$ we have

$$|\rho_l^* - \rho_l|_{k+s} = \left|\int_0^1 \langle \nabla \rho_l(z + tg(z)), g(z)\rangle dt\right|_k \tag{64}$$

$$\leq \sum_{l+m_1+\cdots+m_r = k+s+1} \left|D^l \rho(z+g(z))\right| |D^{m_1} g| \cdots |D^{m_r} g|$$

$$\leq \sum_{m_i \leq k,\ i=1,\ldots,r} \left|D^l \rho(z+g(z))\right| |D^{m_1} g| \cdots |D^{m_r} g| + \sum_{\exists j,\ m_j > k} \left|D^l \rho(z+g(z))\right| |D^{m_1} g| \cdots |D^{m_r} g|$$

$$\leq C(k)\left[(1+|\rho|_{k+s+1})\,|g|_k + (1+|\rho|_{k+1})^{P(k)}\,|g|_{k+s}\right].$$

Rewriting equation (61) for $P',^*$ we obtain

$$[P',^*] = \left[\frac{\partial \rho^*}{\partial \bar{z}}\right]^\tau \cdot \left[I + \left(\left[\frac{\partial \rho^*}{\partial z}\right]\left[\frac{\partial \rho^*}{\partial \bar{z}}\right]^\tau - I\right)\right]^{-1}.$$

Using estimates (63) and (64) we obtain a representation

$$\left[\frac{\partial \rho^*}{\partial z}\right] = \left[\frac{\partial \rho}{\partial z} + B\right],$$

with matrix $B$ satisfying

$$|B|_k \leq C(k)(1 + |\rho|_{k+1})^{P(k)} |g|_k, \qquad (65)$$
$$|B|_{k+s} \leq C(k)\left[(1 + |\rho|_{k+s+1}) |g|_k + (1 + |\rho|_{k+1})^{P(k)} |g|_{k+s}\right].$$

Using this representation in the formula for $[P',^*]$ and estimates (65) we obtain

$$[P',^*] = \left[\frac{\partial \rho}{\partial \bar{z}} + \bar{B}\right]^\tau \cdot [I + A + C]^{-1} \qquad (66)$$

$$= \left[\frac{\partial \rho}{\partial \bar{z}}\right]^\tau \cdot [I + A]^{-1} \cdot [I + A] \cdot [I + A + C]^{-1} + \bar{B}^\tau \cdot [I + A + C]^{-1}$$

$$= \left[\frac{\partial \rho}{\partial \bar{z}}\right]^\tau \cdot [I + A]^{-1} \cdot [I + D] + \bar{B}^\tau \cdot [I + D]$$

$$= [P'] + \left[\frac{\partial \rho}{\partial \bar{z}}\right]^\tau \cdot [I + A]^{-1} \cdot D + \bar{B}^\tau \cdot [I + D],$$

with the matrix $A$ from (59) and matrices

$$C = B \cdot \left[\frac{\partial \rho}{\partial \bar{z}}\right]^\tau + \left[\frac{\partial \rho}{\partial z}\right] \cdot \overline{B}^\tau + B \cdot \overline{B}^\tau, \text{ and } D = [I + A + C]^{-1} - I,$$

satisfying estimates

$$|C|_k, \ |D|_k \leq C(k)(1 + |\rho|_{k+1})^{P(k)} |g|_k, \qquad (67)$$
$$|C|_{k+s}, \ |D|_{k+s} \leq C(k)(1 + |\rho|_{k+1})^{P(k)}\left[(1 + |\rho|_{k+s+1}) |g|_k + |g|_{k+s}\right].$$

Using these estimates in (66) we obtain estimates (62). $\square$

In the proposition below we describe the transformation of a form $\mu \in C^k\left(\mathbf{M}, T'(\mathbf{G})|_{\mathbf{M}}\right)$ under a diffeomorphism $F : \mathbf{M} \to \mathbf{G}$ close to identity. We prove necessary estimates for the new form $\mu^*$, which measures "deviation" of the original deformed CR structure from the inherited CR structure on $F(\mathbf{M})$.

**Proposition 5.4.** *Let $s, k, p \in \mathbb{Z}$ be such that $2 \leq s \leq k$ and $k + s < p$ and let manifold $\mathbf{M} \subset \mathbb{B}(1)$ of the class $C^p$ be defined by functions*

$$\{\rho_l(z) = y_{n-m+l} + r_{n-m+l}(z)\}_1^m$$

*such that $|r_l|_1 < \delta \ll 1$. Let $T''_\mu$ be an almost CR structure on $\mathbb{M}$ defined by the differential form $\mu \in C^p(\mathbb{B}(1))$ such that $|\mu|_k < \epsilon \ll \delta$ and let the functions $\{f_i\}_1^n \in C^p(\mathbb{B}(1))$ be such that $|f|_{k+1} < \epsilon \ll \delta$.*

*Then for small enough $\delta$ and $\epsilon$ there exists a constant $C(k)$ such that for a map*

$$F_i(z) = z_i + f_i(z), \quad i = 1, \ldots, n,$$

and form $\mu^*$ on $\mathbf{M}^* = F(\mathbf{M})$ such that
$$DF\left[T''_\mu(\mathbf{M})\right] = T''_{\mu^*}(\mathbf{M}^*),$$
the following estimates hold
$$|\mu^*|_k \leq C(k)\left(1 + |\rho|_{k+2} + |f|_{k+1} + |\mu|_k\right)^{P(k)} \left(|\bar{\partial}_\mathbf{M} f - \mu|_k + |f|_{k+1} \cdot |\mu|_k\right), \tag{68}$$
$$|\mu^*|_{k+s}$$
$$\leq C(k)\left(1 + |\rho|_{k+2} + |f|_{k+1} + |\mu|_k\right)^{P(k)} \left\{|\bar{\partial}_\mathbf{M} f - \mu|_{k+s} + |f|_{k+s+1} \cdot |\mu|_k + |f|_{k+1} \cdot |\mu|_{k+s}\right.$$
$$\left. + \left(1 + |\rho|_{k+s+2} + |f|_{k+s+1} + |\mu|_{k+s}\right)\left(|\bar{\partial}_\mathbf{M} f - \mu|_k + |f|_{k+1} \cdot |\mu|_k\right)\right\},$$
where $P(k)$ is a polynomial in $k$.

**Proof.** Considering $DF(\bar{X}_i)$ for
$$\bar{X}_i = \bar{Z}_i - \tau \circ \mu(\bar{Z}_i) = \bar{Z}_i + \sum_{j,l} \mu_i^j \frac{\partial \rho_l}{\partial z_j} P''_l - \sum_j \mu_i^j \frac{\partial}{\partial z_j} \in T''_\mu(\mathbf{M}), \tag{69}$$
we obtain
$$DF[\bar{X}_i] = DF\left[\bar{Z}_i + \sum_{j,l} \mu_i^j \frac{\partial \rho_l}{\partial z_j} P''_l - \sum_j \mu_i^j \frac{\partial}{\partial z_j}\right] \tag{70}$$
$$= \sum_s \bar{Z}_i(\bar{F}_s) \frac{\partial}{\partial \bar{z}_s} + \sum_{j,s,l} \mu_i^j \frac{\partial \rho_l}{\partial z_j} P''_l(\bar{F}_s) \frac{\partial}{\partial \bar{z}_s} - \sum_{j,s} \mu_i^j \frac{\partial \bar{F}_s}{\partial z_j} \frac{\partial}{\partial \bar{z}_s}$$
$$+ \sum_s \bar{Z}_i(F_s) \frac{\partial}{\partial z_s} + \sum_{j,s,l} \mu_i^j \frac{\partial \rho_l}{\partial z_j} P''_l(F_s) \frac{\partial}{\partial z_s} - \sum_{j,s} \mu_i^j \frac{\partial F_s}{\partial z_j} \frac{\partial}{\partial z_s}$$
$$= \frac{\partial}{\partial \bar{z}_i} - \sum_s \sum_{l=1}^m \frac{\partial \rho_l}{\partial \bar{z}_i} \bar{p}_l^s \frac{\partial}{\partial \bar{z}_s}$$
$$+ \sum_s \left(\frac{\partial \bar{f}_s}{\partial \bar{z}_i} - \sum_{l=1}^m \frac{\partial \rho_l}{\partial \bar{z}_i} \cdot \sum_{k=1}^n \bar{p}_l^k \frac{\partial \bar{f}_s}{\partial \bar{z}_k} + \sum_{j,l} \mu_i^j \frac{\partial \rho_l}{\partial z_j} P''_l(\bar{F}_s) - \sum_{j,s} \mu_i^j \frac{\partial \bar{f}_s}{\partial z_j}\right) \frac{\partial}{\partial \bar{z}_s}$$
$$+ \sum_s \left(\bar{Z}_i(f_s) - \mu_i^s + \sum_{j,l} \mu_i^j \frac{\partial \rho_l}{\partial z_j} P''_l(f_s) - \sum_j \mu_i^j \frac{\partial f_s}{\partial z_j}\right) \frac{\partial}{\partial z_s}$$
$$= \bar{Z}_i(z) + \beta_i + \sum_s \left(\bar{Z}_i(f_s) - \mu_i^s + \sum_{j,l} \mu_i^j \frac{\partial \rho_l}{\partial z_j} P''_l(f_s) - \sum_j \mu_i^j \frac{\partial f_s}{\partial z_j}\right) \frac{\partial}{\partial z_s},$$
with $\beta_i$ satisfying estimates
$$|\beta_i|_k \leq C(k)\left(1 + |\rho|_{k+1} + |f|_{k+1}\right)^{P(k)} \cdot \left(|\mu|_k + |f|_{k+1}\right), \tag{71}$$
$$|\beta_i|_{k+s} \leq C(k)\left(1 + |\rho|_{k+1} + |f|_{k+1}\right)^{P(k)} \cdot \left[\left(1 + |\rho|_{k+s+1}\right) \cdot \left(|\mu|_k + |f|_{k+1}\right) + |\mu|_{k+s} + |f|_{k+s+1}\right].$$

Tangent $T''(\mathbf{M}^*)$ - vector fields are defined as
$$\bar{Z}_i^*(F(z)) = \frac{\partial}{\partial \bar{z}_i} - \sum_l \frac{\partial \rho_l^*}{\partial \bar{z}_i}(F(z)) \cdot \left(\sum_s \bar{p}_l^{*,s}(F(z)) \frac{\partial}{\partial \bar{z}_s}\right) \tag{72}$$
$$= \frac{\partial}{\partial \bar{z}_i} - \sum_l \frac{\partial \rho_l}{\partial \bar{z}_i}(z) \cdot \left(\sum_s \bar{p}_l^s(z) \frac{\partial}{\partial \bar{z}_s}\right) + \sum_l \left(\frac{\partial \rho_l}{\partial \bar{z}_i}(z) - \frac{\partial \rho_l^*}{\partial \bar{z}_i}(F(z))\right) \cdot \left(\sum_s \bar{p}_l^s(z) \frac{\partial}{\partial \bar{z}_s}\right)$$
$$+ \sum_l \frac{\partial \rho_l^*}{\partial \bar{z}_i}(F(z)) \cdot \left(\sum_s [\bar{p}_l^s(z) - \bar{p}_l^{*,s}(F(z))] \frac{\partial}{\partial \bar{z}_s}\right) = \bar{Z}_i(z) + \nu_i$$

with

$$\nu_i = \sum_l \left(\frac{\partial \rho_l}{\partial \bar{z}_i}(z) - \frac{\partial \rho_l^*}{\partial \bar{z}_i}(F(z))\right) \cdot \left(\sum_s \bar{p}_l^s(z) \frac{\partial}{\partial \bar{z}_s}\right)$$

$$+ \sum_l \frac{\partial \rho_l^*}{\partial \bar{z}_i}(F(z)) \cdot \left(\sum_s \left[\bar{p}_l^s(z) - \bar{p}_l^{*,s}(F(z))\right] \frac{\partial}{\partial \bar{z}_s}\right).$$

Applying Proposition 5.1 and Lemmas 5.2 and 5.3, we obtain estimates for $\nu_i$

$$|\nu_i|_k \le C(k) \left(1 + |\rho|_{k+2} + |f|_{k+1}\right)^{P(k)} \cdot |f|_{k+1}, \tag{73}$$

$$|\nu_i|_{k+s} \le C(k) \left(1 + |\rho|_{k+2} + |f|_{k+1}\right)^{P(k)} \left[(1 + |\rho|_{k+s+2}) \cdot |f|_{k+1} + \cdot |f|_{k+s+1}\right].$$

From the choice of functions $\{\rho_l\}_1^m$ and estimates (73) we conclude that vectors $\{\bar{Z}_i\}_1^{n-m}$ represent a basis in $T''(\mathbf{M})$ and for $\epsilon$ small enough vectors $\{\bar{Z}_i^*\}_1^{n-m}$ represent a basis in $T''(\mathbf{M}^*)$.

Combining (70) and (72) we obtain that for $\epsilon$ small enough

$$T''_{\mu^*}(\mathbf{M}^*) = DF\left[T''_\mu(\mathbf{M})\right] = \mathrm{Span}\bigg\{\bar{Z}_i^* + \sum_s \lambda_i^s \frac{\partial}{\partial \bar{z}_s} \tag{74}$$

$$+ \sum_s \left(\bar{Z}_i(f_s) - \mu_i^s + \sum_{j,l} \mu_i^j \frac{\partial \rho_l}{\partial z_j} P_l''(f_s) - \sum_j \mu_i^j \frac{\partial f_s}{\partial z_j}\right) \frac{\partial}{\partial z_s}\bigg\},$$

$$= \mathrm{Span}\bigg\{\bar{Z}_i^* + \sum_s \lambda_i^s \frac{\partial}{\partial \bar{z}_s} + \sum_s \left(\bar{Z}_i(f_s) - \mu_i^s + \sum_{j,l} \mu_i^j \frac{\partial \rho_l}{\partial z_j} P_l''(f_s) - \sum_j \mu_i^j \frac{\partial f_s}{\partial z_j}\right) \frac{\partial}{\partial z_s}\bigg\}_1^{n-m},$$

with $\lambda_i$ satisfying estimates

$$|\lambda_i|_k \le C(k) \left(1 + |\rho|_{k+2} + |f|_{k+1}\right)^{P(k)} \cdot \left(|\mu|_k + |f|_{k+1}\right), \tag{75}$$

$$|\lambda_i|_{k+s} \le C(k) \left(1 + |\rho|_{k+2} + |f|_{k+1}\right)^{P(k)} \cdot \left[(1 + |\rho|_{k+s+2}) \cdot (|\mu|_k + |f|_{k+1}) + |\mu|_{k+s} + |f|_{k+s+1}\right].$$

Using formulas

$$\begin{bmatrix} \bar{Z}_1^* \\ \vdots \\ \bar{Z}_{n-m}^* \\ P_1^{*,''} \\ \vdots \\ P_m^{*,''} \end{bmatrix} = \begin{bmatrix} I & -\{\frac{\partial \rho_i}{\partial \bar{z}_j}\}_{j=1,\ldots,n-m}^{i=1,\ldots,m} \\ 0 & I \end{bmatrix} \cdot \begin{bmatrix} \frac{\partial}{\partial \bar{z}_1} & \cdots & 0 & 0 & \cdots & 0 \\ \vdots & \vdots & \vdots & \vdots & \vdots & \vdots \\ 0 & \cdots & \frac{\partial}{\partial \bar{z}_{n-m}} & 0 & \cdots & 0 \\ \bar{p}_1^1 \frac{\partial}{\partial \bar{z}_1} & \cdots & \cdots & \cdots & \cdots & \bar{p}_1^n \frac{\partial}{\partial \bar{z}_n} \\ \vdots & \vdots & \vdots & \vdots & \vdots & \vdots \\ \bar{p}_m^1 \frac{\partial}{\partial \bar{z}_1} & \cdots & \cdots & \cdots & \cdots & \bar{p}_m^n \frac{\partial}{\partial \bar{z}_n} \end{bmatrix}$$

and

$$\begin{bmatrix} \frac{\partial}{\partial \bar{z}_1} & \cdots & 0 & 0 & \cdots & 0 \\ \vdots & \vdots & \vdots & \vdots & \vdots & \vdots \\ 0 & \cdots & \frac{\partial}{\partial \bar{z}_{n-m}} & 0 & \cdots & 0 \\ \bar{p}_1^1 \frac{\partial}{\partial \bar{z}_1} & \cdots & \cdots & \cdots & \cdots & \bar{p}_1^n \frac{\partial}{\partial \bar{z}_n} \\ \vdots & \vdots & \vdots & \vdots & \vdots & \vdots \\ \bar{p}_m^1 \frac{\partial}{\partial \bar{z}_1} & \cdots & \cdots & \cdots & \cdots & \bar{p}_m^n \frac{\partial}{\partial \bar{z}_n} \end{bmatrix}$$

$$= \begin{bmatrix} 1 & \cdots & 0 & 0 & \cdots & 0 \\ \vdots & \vdots & \vdots & \vdots & \vdots & \vdots \\ 0 & \cdots & 1 & 0 & \cdots & 0 \\ \bar{p}_1^1 & \cdots & \cdots & \cdots & \cdots & \bar{p}_1^n \\ \vdots & \vdots & \vdots & \vdots & \vdots & \vdots \\ \bar{p}_m^1 & \cdots & \cdots & \cdots & \cdots & \bar{p}_m^n \end{bmatrix} \cdot \begin{bmatrix} \frac{\partial}{\partial \bar{z}_1} & \cdots & 0 \\ \vdots & \vdots & \vdots \\ 0 & \cdots & \frac{\partial}{\partial \bar{z}_n} \end{bmatrix},$$

and Lemmas 5.2 and 5.3 we obtain formula

$$\begin{bmatrix} \frac{\partial}{\partial \bar{z}_1} \\ \vdots \\ \frac{\partial}{\partial \bar{z}_n} \end{bmatrix} = [I + A] \cdot \begin{bmatrix} \bar{Z}_i^* \\ P^{*,\prime\prime} \end{bmatrix} \quad (76)$$

with a matrix $A$ satisfying estimates

$$|A|_0 \leq C\delta, \quad (77)$$

$$|A|_k \leq C(k) \left(1 + |\rho|_{k+1}\right)^{P(k)},$$

$$|A|_{k+s} \leq C(k) \left(1 + |\rho|_{k+1}\right)^{P(k)} \left(1 + |\rho|_{k+s+1}\right).$$

Using formula (76) and estimates (75) and (77), we obtain equality

$$\begin{bmatrix} \bar{Z}_i^* + \sum_s \lambda_i^s \frac{\partial}{\partial \bar{z}_s} \\ P^{*,\prime\prime} \end{bmatrix} = \begin{bmatrix} I + B & C \\ 0 & I \end{bmatrix} \cdot \begin{bmatrix} \bar{Z}_i^* \\ P^{*,\prime\prime} \end{bmatrix}$$

with matrices $B, C$ satisfying

$$|B|_k, |C|_k \leq C(k) \left(1 + |\rho|_{k+2} + |f|_{k+1}\right)^{P(k)} \cdot \left(|\mu|_k + |f|_{k+1}\right), \quad (78)$$

$$|B|_{k+s}, |C|_{k+s}$$

$$\leq C(k) \left(1 + |\rho|_{k+2} + |f|_{k+1}\right)^{P(k)} \cdot \left[\left(1 + |\rho|_{k+s+2}\right) \cdot \left(|\mu|_k + |f|_{k+1}\right) + |\mu|_{k+s} + |f|_{k+s+1}\right].$$

Comparing the last equality with

$$\begin{bmatrix} \bar{Z}_i^* + \sum_{s,l} \mu_i^{*,s} \frac{\partial \rho_l^*}{\partial z_s} P_l^{*,\prime\prime} \\ P^{*,\prime\prime} \end{bmatrix} = \begin{bmatrix} I & D \\ 0 & I \end{bmatrix} \cdot \begin{bmatrix} \bar{Z}_i^* \\ P^{*,\prime\prime} \end{bmatrix}$$

we obtain

$$\begin{bmatrix} \bar{Z}_i^* + \sum_{s,l} \mu_i^{*,s} \frac{\partial \rho_l}{\partial z_s} P_l^{*,\prime\prime} \\ P^{*,\prime\prime} \end{bmatrix} = \begin{bmatrix} [I+B]^{-1} & -[I+B]^{-1}C + D \\ 0 & I \end{bmatrix} \cdot \begin{bmatrix} \bar{Z}_i^* + \sum_s \lambda_i^s \frac{\partial}{\partial \bar{z}_s} \\ P^{*,\prime\prime} \end{bmatrix}. \quad (79)$$

Using then formula (69) for $\bar{X}_i^*$:

$$\bar{X}_i^* = \bar{Z}_i^* + \sum_{s,l} \mu_i^{*,s} \frac{\partial \rho_l^*}{\partial z_s} P_l^{*,\prime\prime} - \sum_s \mu_i^{*,s} \frac{\partial}{\partial z_s}$$

and equality (79) we obtain formula

$$[\mu_i^{*,s}] = -[I+B]^{-1} \cdot \left[\bar{Z}_i(f_s) - \mu_i^s + \sum_{j,l} \mu_i^j \frac{\partial \rho_l}{\partial z_j} P_l^{\prime\prime}(f_s) - \sum_j \mu_i^j \frac{\partial f_s}{\partial z_j}\right].$$

Using in the last formula estimates (78) we obtain estimates (68). □

## 6. Proof of Theorem 1.

Let now $\mathcal{E}_0(\mathbb{M}) = \mathbf{M}_0 \subset \mathbf{G}$ be a fixed compact, regular 3-pseudoconcave CR submanifold such that $\dim H^1\left(\mathbf{M}_0, T'(\mathbf{G})|_{\mathbf{M}_0}\right) = 0$. In order to prove Theorem 1 we have to construct for any sufficiently small deformation $T''_\mu(\mathbf{M}_0)$, defined by a differential form $\mu \in C^p_{(0,1)}\left(\mathbf{M}_0, T'(\mathbf{G})|_{\mathbf{M}_0}\right)$ an appropriate embedding $\mathcal{E} : \mathbb{M} \to \mathbf{G}$ such that $T''_\mu(\mathbf{M}_0) = \mathcal{E}^*\left(T''(\mathbf{G})\right)$. We construct map $\mathcal{E}$ in a form $\mathcal{E} = \mathcal{F} \circ \mathcal{E}_0$ with $\mathcal{F}$ being a limit of a sequence of maps $\{\mathcal{F}^{(j)} : \mathbf{M}_0 \to \mathbf{G}\}_0^\infty$, close to identity.

We consider parameterization of the set of $C^l$-diffeomorphisms $\mathrm{Diff}^l(\mathbf{M}_0, \mathbf{G})$ close to identity given by the exponential map [S], [L]:

$$\exp : C^l_{(0,0)}\left(\mathbf{M}_0, T'(\mathbf{G})|_{\mathbf{M}_0}\right) \to \mathrm{Diff}^l(\mathbf{M}_0, \mathbf{G}).$$

If $\{\mathcal{U}_j\}_1^N$ is a finite covering of the neighborhood of submanifold $\mathbf{M}_0 \subset \mathbf{G}$ by coordinate neighborhoods with coordinates $\{z_1^j, \ldots, z_n^j\}$, then for a vector field $\xi \in C^l_{(0,0)}\left(\mathbf{M}_0, T'(\mathbf{G})|_{\mathbf{M}_0}\right)$ with $|\xi|_l$ small enough and a local representation $\xi = (\xi_1, \ldots, \xi_n)$ in $\mathcal{U} \in \{\mathcal{U}_j\}_1^N$ we have

$$[\exp \xi(z)]_i = z_i + \xi_i(z) + A_i(z, \xi(z)), \tag{80}$$

with

$$\begin{aligned} |A|_{k+2} &< C|\xi|_{k+2}^2, \\ |A|_{k+s+2} &< C|\xi|_{k+s+2}|\xi|_{k+2}, \end{aligned} \tag{81}$$

for $0 \leq s \leq k - 4$ and some $C > 0$.

To construct sequence $\{\mathcal{F}^{(j)}\}_0^\infty$ we:
(1) fix a family of smoothing operators ([W2],[M]) for $0 < t < \infty$

$$S_t : C^l_{(0,i)}(\mathbf{G}, T'(\mathbf{G})) \to C^\infty_{(0,i)}(\mathbf{G}, T'(\mathbf{G}))$$

such that for $0 \leq q \leq r \leq p$

$$\begin{aligned} |S_t f|_r &\leq C t^{q-r} |f|_q, \\ |(I - S_t) f|_q &\leq C t^{r-q} |f|_r, \end{aligned} \tag{82}$$

with some $C > 0$,
(2) consider a sequence of numbers

$$t_{j+1} = (t_j)^{7/6},$$

(3) consider the following sequence of maps $\{F^{(j)}\}_0^\infty$
   (i) $F^{(j)} = \exp\left(\xi^{(j)}\right)$
   (ii) $\xi^{(j)} = S_{t_j} \circ P^{(j)}(\mu^{(j)})$,
   (iii) $\mu^{(j+1)} = \left(\mu^{(j)}\right)^* \in C^p_{(0,1)}\left(\mathbf{M}_{j+1}, T'(\mathbf{G})|_{\mathbf{M}_{j+1}}\right)$,
   where $\mathbf{M}_{j+1} = F^{(j)}(\mathbf{M}_j)$, $P^{(j)} = P_{\mathbf{M}_j}$ from Theorem 2 and the construction of $\left(\mu^{(j)}\right)^*$ for given $\mu^{(j)}$ and $F^{(j)}$ is described in Proposition 5.4.

Then we define sequence $\mathcal{F}^{(j)}$ as

$$\mathcal{F}^{(j)} = F^{(j)} \circ F^{(j-1)} \circ \cdots \circ F^{(0)}$$

and obtain Theorem 1 as a corollary of the following

**Proposition 6.1.** *Let $k, p \in \mathbb{Z}$ be such that $17 \le k$, $2k \le p - 5$ and let $\mathbf{M}_0 \subset \mathbf{G}$ be a regular 3-pseudoconcave submanifold of the class $C^p$. Then there exist $t_0, \epsilon > 0$ such that for $\mu^{(0)} \in C^p_{(0,1)}\left(\mathbf{M}_0, T'(\mathbf{G})|_{\mathbf{M}_0}\right)$, satisfying integrability condition and such that $\delta_0(2k) := |\mu^{(0)}|_{2k} < \epsilon$ the following conditions are satisfied for $\theta_j = |\xi^{(j)}|_{k+2}$ and $\delta_j = |\mu^{(j)}|_k$*

(a) $\sum_0^\infty \theta_j < \infty$,
(b) $\lim_{j \to \infty} \delta_j = 0$.

**Proof.** We obtain the statement of Proposition 6.1 as a corollary of the following estimates
$$1 + |\rho^{(j)}|_{k+2} \le C,$$
$$t_j^{-6} \delta_j \le 1/2, \qquad (83)$$
$$t_j^{s-12}\left[|\rho^{(j)}|_{k+s+2} + \delta_j(s)\right] \le 1/4,$$
where $13 \le s \le k - 4$, $\delta_j(s) = |\mu^{(j)}|_{k+s}$, $\rho^{(j)} = \rho^{(j-1)} \circ \mathcal{G}^{(j)}$ and $\mathcal{G}^{(j)} = \left(\mathcal{F}^{(j)}\right)^{-1}$.

Statement (a) of Proposition 6.1 follows from the first two estimates in (83) by the following sequence of inequalities:
$$t_j^{-1}\theta_j = t_j^{-1}|S_{t_j} \circ P^{(j)}(\mu^{(j)})|_{k+2} \le C t_j^{-6}|P^{(j)}(\mu^{(j)})|_{k-3} \qquad (84)$$
$$\le C t_j^{-6} \delta_j \left(1 + |\rho^{(j)}|_{k+2}\right)^{P(k)} \le C t_j^{-6} \delta_j \le C,$$
where we used also estimates (82) and the first estimate in (3).

Statement (b) obviously follows from the second estimate in (83).

The choice of $t_0$ will be defined in the lemma below. The choice of $\epsilon = \frac{1}{2} t_0^6$ guarantees the (83) for the initial step of induction.

We prove estimates (83) by induction, assuming that (83) holds for some fixed $j$ and prove it for $j + 1$. In the proof of the step of induction we will need the following lemma.

**Lemma 6.2.** *In the notation of Proposition 6.1 and assuming that estimates (83) hold for $j$ the following estimates hold for $j + 1$*

(i) $\delta_{j+1} \le C \left(1 + |\rho^{(j)}|_{k+2}\right)^{P(k)} \left\{t_j^{-4}\delta_j^2 + t_j^{s-4}\left[\left(1 + |\rho^{(j)}|_{k+s+2}\right)\delta_j + \delta_j(s)\right]\right\}$,

(ii) $|\rho^{(j+1)}|_{k+s+2} \le C \left(1 + |\rho^{(j)}|_{k+2}\right)^{P(k)} \left[|\rho^{(j)}|_{k+s+2} + t_j^{-s-5}\delta_j\right]$,

(iii) $\delta_{j+1}(s) \le C \left(1 + |\rho^{(j)}|_{k+2}\right)^{P(k)} \left\{\delta_j(s) + t_j^{-s-4}\delta_j + t_j^{-4}\delta_j \left|\rho^{(j)}\right|_{k+s+2}\right\}$.

**Proof.** Transforming representation (2) into
$$\mu = \bar{\partial}_{\mathbf{M}} \circ S_t \circ P(\mu) + \bar{\partial}_{\mathbf{M}} \circ (I - S_t) \circ P(\mu) + S_t \circ Q \circ \bar{\partial}_{\mathbf{M}}(\mu) + (I - S_t) \circ Q \circ \bar{\partial}_{\mathbf{M}}(\mu) \qquad (85)$$
and using estimates
$$\left|\bar{\partial}_{\mathbf{M}} \circ \left(I - S_{t_j}\right) \circ P^{(j)}(\mu^{(j)})\right|_k \le \left|\left(I - S_{t_j}\right) \circ P^{(j)}(\mu^{(j)})\right|_{k+1} \le C t_j^{s-4}\left|P^{(j)}(\mu^{(j)})\right|_{k+s-3}$$
$$\le C t_j^{s-4}\left(1 + |\rho^{(j)}|_{k+2}\right)^{P(k)}\left[\left(1 + |\rho^{(j)}|_{k+s+2}\right)\left|\mu^{(j)}\right|_k + \left|\mu^{(j)}\right|_{k+s}\right],$$
$$\left|S_{t_j} \circ Q^{(j)} \circ \bar{\partial}_{\mathbf{M}}(\mu^{(j)})\right|_k \le C t_j^{-4}\left|Q^{(j)} \circ \bar{\partial}_{\mathbf{M}}(\mu^{(j)})\right|_{k-4}$$
$$\le C t_j^{-4}\left(1 + |\rho^{(j)}|_{k+1}\right)^{P(k)}\left|\bar{\partial}_{\mathbf{M}}(\mu^{(j)})\right|_{k-1},$$
and
$$\left|\left(I - S_{t_j}\right) \circ Q^{(j)} \circ \bar{\partial}_{\mathbf{M}}(\mu^{(j)})\right|_k \le C t_j^{s-4}\left|Q^{(j)} \circ \bar{\partial}_{\mathbf{M}}(\mu^{(j)})\right|_{k+s-4}$$

$$\leq C t_j^{s-4} \left(1 + |\rho^{(j)}|_{k+1}\right)^{P(k)} \left[\left(1 + |\rho^{(j)}|_{k+s+1}\right) \left|\bar{\partial}_{\mathbf{M}}(\mu^{(j)})\right|_{k-1} + \left|\bar{\partial}_{\mathbf{M}}(\mu^{(j)})\right|_{k+s-1}\right],$$

we obtain
$$\left|\mu^{(j)} - \bar{\partial}_{\mathbf{M}} \circ S_{t_j} \circ P^{(j)}(\mu^{(j)})\right|_k \tag{86}$$
$$\leq C \left(1 + |\rho^{(j)}|_{k+2}\right)^{P(k)} \left\{ t_j^{-4} \left|\bar{\partial}_{\mathbf{M}}(\mu^{(j)})\right|_{k-1} + t_j^{s-4} \left[\left(1 + |\rho^{(j)}|_{k+s+2}\right) \delta_j + \delta_j(s)\right] \right\}.$$

To estimate the first term in the right hand side of (86) we use integrability condition (50) and obtain

$$\left|\bar{\partial}_{\mathbf{M}}(\mu^{(j)})\right|_{k-1} \leq \left[ \left|\sum_{\gamma,\beta} \left(\frac{\partial^\tau \mu_\beta^{(j),\alpha}}{\partial z_\gamma}\right) \mu^{(j),\gamma} \wedge \bar{Z}^\beta \right|_{k-1} \right.$$

$$\left. + \left|\sum_\beta \mu_\beta^{(j),\alpha} \cdot \sum_{l=1}^m \left(\bar{\partial}_{\mathbf{M}} \bar{p}_k^\beta - \sum_\gamma \frac{\partial^\tau \bar{p}_k^\beta}{\partial z_\gamma} \mu^{(j),\gamma}\right) \wedge \left(\sum_\delta \frac{\partial \rho_l}{\partial z_\delta} \mu^{(j),\delta}\right) \right|_{k-1} \right]$$

$$\leq C \left(1 + |\rho^{(j)}|_{k+2}\right)^{P(k)} \left|\mu^{(j)}\right|_k^2.$$

Using this estimate in (86) we obtain
$$\left|\mu^{(j)} - \bar{\partial}_{\mathbf{M}} \xi^{(j)}\right|_k = \left|\mu^{(j)} - \bar{\partial}_{\mathbf{M}} \circ S_{t_j} \circ P^{(j)}(\mu^{(j)})\right|_k \tag{87}$$
$$\leq C \left(1 + |\rho^{(j)}|_{k+2}\right)^{P(k)} \left\{ t_j^{-4} \delta_j^2 + t_j^{s-4} \left[\left(1 + |\rho^{(j)}|_{k+s+2}\right) \delta_j + \delta_j(s)\right] \right\}.$$

Applying then Proposition 5.4 with $f^{(j)}(z) = \xi^{(j)}(z) + A\left(z, \xi^{(j)}(z)\right)$ from (80) and estimates (81), we obtain

$$\delta_{j+1} = \left|\mu^{(j+1)}\right|_k \leq C \left(1 + |\rho^{(j)}|_{k+2} + \delta_j + |f^{(j)}|_{k+1}\right)^{P(k)} \left(\left|\bar{\partial}_{\mathbf{M}} f^{(j)} - \mu^{(j)}\right|_k + \left|\mu^{(j)}\right|_k |f^{(j)}|_{k+1}\right)$$

$$\leq C \left(1 + |\rho^{(j)}|_{k+2} + \delta_j + \left(|\xi^{(j)}|_{k+1} + |\xi^{(j)}|_{k+1}^2\right)\right)^{P(k)}$$
$$\left[\left|\bar{\partial}_{\mathbf{M}} \xi^{(j)} - \mu^{(j)}\right|_k + |\xi^{(j)}|_{k+1}^2 + |\mu^{(j)}|_k \left(|\xi^{(j)}|_{k+1} + |\xi^{(j)}|_{k+1}^2\right)\right]$$
$$\leq C \left(1 + |\rho^{(j)}|_{k+2}\right)^{P(k)} \left\{ t_j^{-4} \delta_j^2 + t_j^{s-4} \left[\left(1 + |\rho^{(j)}|_{k+s+2}\right) \delta_j + \delta_j(s)\right] \right\},$$

where in the last inequality we used estimate (87), estimate
$$\left|\xi^{(j)}\right|_{k+1} = \left|S_{t_j} \circ P^{(j)}(\mu^{(j)})\right|_{k+1} \leq C \left(1 + |\rho^{(j)}|_{k+2}\right)^{P(k)} t_j^{-4} \delta_j, \tag{88}$$

and estimates (83) for $j$.

For part (ii) of the lemma, using estimates (56), (84) and (88) we obtain
$$|\rho^{(j+1)}|_{k+s+2} = \left|\rho^{(j)}\left(z + g^{(j)}(z)\right)\right|_{k+s+2}$$
$$\leq C \left[|\rho^{(j)}|_{k+s+2} \left(1 + |g^{(j)}|_{k+2}\right)^{P(k)} + \left(1 + |\rho^{(j)}|_{k+2}\right)^{P(k)} |g^{(j)}|_{k+s+2}\right]$$
$$\leq C \left(1 + |\rho^{(j)}|_{k+2}\right)^{P(k)} \left[|\rho^{(j)}|_{k+s+2} \left(1 + |f^{(j)}|_{k+2}\right)^{P(k)} + |f^{(j)}|_{k+s+2}\right]$$
$$\leq C \left(1 + |\rho^{(j)}|_{k+2}\right)^{P(k)} \left[|\rho^{(j)}|_{k+s+2} + |\xi^{(j)}|_{k+s+2}\right]$$
$$\leq C \left(1 + |\rho^{(j)}|_{k+2}\right)^{P(k)} \left[|\rho^{(j)}|_{k+s+2} + t_j^{-s-5} \delta_j\right],$$

where in the last inequality we used estimate
$$\left|\xi^{(j)}\right|_{k+s+2} = \left|S_{t_j} \circ P^{(j)}(\mu^{(j)})\right|_{k+s+2} \leq C t_j^{-s-5} \left|P^{(j)}(\mu^{(j)})\right|_{k-3} \tag{89}$$

$$\leq C \left(1 + |\rho^{(j)}|_{k+2}\right)^{P(k)} t_j^{-s-5} \delta_j.$$

To prove part (iii) of the lemma we use second estimate in (68) from Proposition 5.4, estimates (81) and inductive assumptions (83), and obtain

$$\delta_{j+1}(s) = \left|\mu^{(j+1)}\right|_{k+s} \tag{90}$$

$$\leq C \left(1 + |\rho^{(j)}|_{k+2}\right)^{P(k)} \left\{\left|\mu^{(j)}\right|_{k+s} + \left|\bar{\partial}_{\mathbf{M}} f^{(j)}\right|_{k+s} + \left|f^{(j)}\right|_{k+s+1} \cdot \left|\mu^{(j)}\right|_k + \left|f^{(j)}\right|_{k+1} \cdot \left|\mu^{(j)}\right|_{k+s}\right.$$

$$+ \left(1 + \left|\rho^{(j)}\right|_{k+s+2} + \left|f^{(j)}\right|_{k+s+1} + |\mu^{(j)}|_{k+s}\right) \left(\left|\bar{\partial}_{\mathbf{M}} f^{(j)} - \mu^{(j)}\right|_k + \left|f^{(j)}\right|_{k+1} \cdot \left|\mu^{(j)}\right|_k\right)\right\},$$

$$\leq C \left(1 + |\rho^{(j)}|_{k+2}\right)^{P(k)} \left\{\delta_j(s) + t_j^{-s-4}\delta_j + t_j^{-s-4}\delta_j^2 + t_j^{-4}\delta_j \cdot \delta_j(s)\right.$$

$$+ \left(1 + \left|\rho^{(j)}\right|_{k+s+2} + t_j^{-s-4}\delta_j + \delta_j(s)\right) \left(t_j^{-4}\delta_j + \delta_j + t_j^{-4}\delta_j^2\right)\right\},$$

$$\leq C \left(1 + |\rho^{(j)}|_{k+2}\right)^{P(k)} \left\{\delta_j(s) + t_j^{-s-4}\delta_j + t_j^{-4}\delta_j \left|\rho^{(j)}\right|_{k+s+2}\right\}.$$

$\square$

For the $(j+1)$-th induction step in the first estimate in (83) we have

$$|\rho^{(j+1)}(z)|_{k+2} = |\rho^{(j)}(z + g(z))|_{k+2} \leq |\rho^{(j)}|_{k+2} \left(1 + |g^{(j)}|_{k+2}\right)^{P(k)} \leq |\rho^{(j)}|_{k+2} (1 + C\theta_j)^{P(k)}.$$

Combining this estimate with (84) we obtain

$$\left(1 + |\rho^{(j+1)}|_{k+2}\right) \leq \left(1 + |\rho^{(0)}|_{k+2}\right) \prod_1^j (1 + Ct_j)^{P(k)} \leq C.$$

For the second $(j+1)$-th estimate in (83) we use estimates (83) for $j$ and estimate $(i)$ from lemma 6.2 and obtain

$$t_{j+1}^{-6}\delta_{j+1} = t_j^{-7}\delta_{j+1}$$

$$\leq Ct_j^{-7} \left(1 + |\rho^{(j)}|_{k+2}\right)^{P(k)} \left\{t_j^{-4}\delta_j^2 + t_j^{s-4}\left[\left(1 + |\rho^{(j)}|_{k+s+2}\right)\delta_j + \delta_j(s)\right]\right\}$$

$$\leq Ct_j \left\{t_j^{-12}\delta_j^2 + t_j^{s-12}\left[\left(1 + |\rho^{(j)}|_{k+s+2}\right)\delta_j + \delta_j(s)\right]\right\}$$

$$\leq Ct_j \left(t_j^{-6}\delta_j\right)^2 + Ct_j t_j^{s-12}\left[\left(1 + |\rho^{(j)}|_{k+s+2}\right)\delta_j + \delta_j(s)\right] \leq 1/4 + 1/4 = 1/2,$$

where we used conditions $t_j \leq t_0 \leq C^{-1}$.

For the third $(j+1)$-th estimate in (83) we obtain using estimate $(ii)$ from lemma 6.2, second estimate from (83) and estimate (90)

$$t_{j+1}^{s-12} \left[|\rho^{(j+1)}|_{k+s+2} + \delta_{j+1}(s)\right] \leq C \left(1 + |\rho^{(j)}|_{k+2}\right)^{P(k)}$$

$$\times t_j^{\frac{7}{6}(s-12)} \left[|\rho^{(j)}|_{k+s+2} + t_j^{-s-5}\delta_j + \delta_j(s) + t_j^{-s-4}\delta_j + t_j^{-4}\delta_j \left|\rho^{(j)}\right|_{k+s+2}\right]$$

$$\leq Ct_j^{\frac{1}{6}(s-12)} t_j^{s-12} \left[|\rho^{(j)}|_{k+s+2} + \delta_j(s)\right] + Ct_j^{\frac{7}{6}(s-12)} t_j^{-s-5}\delta_j$$

$$\leq \frac{1}{4} Ct_j^{\frac{1}{6}(s-12)} + Ct_j^{\frac{1}{6}(s-12)} t_j^{-5}\delta_j,$$

with expression on the right that can be made arbitrarily small if $s \geq 13$. $\square$

## References

[AiH]  R.A. Airapetian, G.M. Henkin, Integral representation of differential forms on Cauchy-Riemann manifolds and the theory of CR functions, I, II, Russ. Math. Surv., 39 (1984), 41-118, Math.USSR Sbornik 55:1 (1986), 91-111.

[B]    L. Boutet de Monvel, Intgration des quations de Cauchy-Riemann induites formelles, Sminaire Goulaouic-Lions-Schwartz 1974-1975, Exp. No. 9.

[BT]   M.S. Baouendi, F. Treves, A property of the functions and distributions annihilated by a locally integrable system of complex vector fields, Ann. of Math. 113 (1981), 341-421.

[Ca]   D. Catlin, Sufficient conditions for the extension of CR structures. J. Geom. Anal. 4 (1994), no. 4,467-538.

[EH]   C. Epstein, G. Henkin, Stability of embeddings for pseudoconcave surfaces and their boundaries, Acta Math. 185 (2000), 161-237.

[Ha1]  R.S. Hamilton, The inverse function theorem of Nash and Moser, Bull.Amer.Math.Soc. 7:1 (1982), 65-222.

[Ha2]  R.S. Hamilton, Deformation of complex structures on manifolds with boundary, I J. Diff. Geometry 12 (1977), 1-45; II J.Diff. Geometry 14 (1979), 409-473; III Preprint.

[He1]  G.M. Henkin, The Hans Lewy equation and analysis on pseudoconvex manifolds, Math.USSR Sbornik 31 (1977), 59-130.

[He2]  G.M. Henkin, Solution des equations de Cauchy-Riemann tangentielle sur les varietes de Cauchy-Riemann q-concaves, C. R. Acad. Sci. Paris, 292 (1981), 27-30.

[Ki]   G.K. Kiremidjian, Extendible Pseudocomplex Structures I, J. of Approximation Theory 19 (1977), 281-303; II J. D'Analyse Mathématique 30 (1976), 304-329.

[KN]   J.J. Kohn, L. Nirenberg, Non-coercive boundary value problems, Comm. Pure Appl. Math. 18 (1965), 443-492.

[Ku]   M. Kuranishi, Application of $\bar{\partial}_b$ to deformation of isolated singularities, Proc. Symp. Pure Math. 30:1 (1977), 97-103.

[L]    S. Lang, Differential manifolds, Springer-Verlag, 1985.

[M]    J.K. Moser, A rapidly convergent iteration method and nonlinear differential equations II, Ann. Scuola Norm. Sup. Pisa 20 (1966), 499-535.

[P1]   P. Polyakov, Banach cohomology on piecewise strictly pseudoconvex domains, Mat. Sb. (N.S.) 88 (130) (1972), 238-255.

[P2]   P. Polyakov, Sharp estimates for operator $\bar{\partial}_{\mathbf{M}}$ on a q-concave CR manifold, The Journal of Geometric Analysis 6:2 (1996), 233-276.

[P3]   P. Polyakov, Sharp Lipschitz estimates for operator $\bar{\partial}_{\mathbf{M}}$ on a q-concave CR manifold, Preprint, http://xxx.lanl.gov/math.CV/0012250.

[P4]   P. Polyakov, Global $\bar{\partial}_{\mathbf{M}}$-homotopy with $C^k$ estimates for a family of compact, regular q-pseudoconcave CR manifolds, Preprint, http://xxx.lanl.gov/math.CV/0109180.

[RS]   M. Range, Y-T. Siu, Uniform estimates for the $\bar{\partial}$-equation on domains with piecewise smooth strictly pseudoconvex boundaries. Math. Ann. 206 (1973), 325-354.

[Siu]  Y.-T. Siu, The $\bar{\partial}$ problem with uniform bounds on derivatives, Math.Ann. 207 (1974), 163-176.

[S]    S. Sternberg, Lectures on differential geometry, Chelsea, N.Y., 1983.

[W1]   S. Webster, On the local solution of the tangential Cauchy-Riemann equations, Ann.Inst.Henri Poincare, 6:3, (1989),167-182.

[W2]   S. Webster, On the proof of Kuranishi's embedding theorem, Ann.Inst.Henri Poincare, 6:3, (1989), 183-207.

Department of Mathematics, University of Wyoming, Laramie, WY 82071
*E-mail address*: polyakov@uwyo.edu